\input amstex.tex
\input amsppt.sty
\TagsAsMath

\magnification=1200
\hsize=5.0in\vsize=7.in
\hoffset=0.2in\voffset=0cm
\nonstopmode

\def\th{\text{th}}
\document

\topmatter

\title{On asymptotic stability in 3D of kinks for the $\phi ^4$ model}
\endtitle

\author
Scipio Cuccagna
\endauthor

\address
Scipio Cuccagna, Dipartimento di Scienze e  Metodi per l'Ingegneria,
Universita di Modena e Reggio Emilia, via Fogliani 1, Reggio Emilia
42100, Italy
\endaddress

\email cuccagna.scipio\@unimore.it \endemail

\abstract  We add to a kink, which is a 1 dimensional structure,
two transversal directions. We then check its  asymptotic
stability with respect to compactly supported perturbations in 3D
and a time evolution under a Nonlinear Wave Equation (NLW). The
problem is inspired by work by Jack Xin  on asymptotic stability
in dimension larger than 1 of fronts for reaction diffusion
equations. The proof involves a separation of variables. The
transversal variables are treated as in work on Nonlinear Klein
Gordon Equation (NLKG) originating from
  Klainerman and from Shatah in a particular elaboration
due to Delort {\it{et al}}.
The longitudinal variable is treated by means of a result by Weder
on dispersion for Schr\"odinger operators in 1D.
\endabstract

\thanks Research fully supported by a special grant from the Italian
Ministry of Education, University and Research.
\endthanks
\endtopmatter

\head \S 1 Introduction \endhead
\bigskip

\noindent Set $\vec{x} =(x,y) \in \Bbb R \times \Bbb R^2 $,
$\Delta $ the full Laplacian, $\Delta_y$ the Laplacian in the
variables $y$. We consider the NLW
$$u_{tt} -\Delta u-u+u^3=0\,,\,(t,\vec{x})
\in \Bbb R \times \Bbb R^n \,,\, n=3\tag 1.1$$
and the kink solution
$\th (x)=\tanh (2^{-\frac 12}x).$
We consider solutions of the form
$$u(t,x,y)= \th (x)+w(t,x,y),\tag 1.2$$
with initial data (chosen real valued
not to complicate notation)
$$w(0,x,y)= w_0(x,y)
\quad,\, w_t(0,x,y)= w_1(x,y).\tag 1.3$$
We prove:

\proclaim{Theorem 1.1} Consider equation (1.1) with solution (1.2)
with initial data (1.3). Assume that the $ w_j(x,y)$ are smooth,
supported in $|x|^2+|y|^2\le K^2$, for $K>0$ some fixed constant.
Then there exists $\epsilon_0$ such that if $s=24$ and if $( w_0,
w_1)\in H^s(\Bbb R^3) \times H^{s-1}(\Bbb R^3)$ with norm smaller
than $\epsilon_0$, then the solution of the form (1.2) is globally
defined in time and such that
 $$| w (t,\vec{x} )|\sqrt {1+t}\sqrt {1+|t-|y ||}
 < \infty .$$
\endproclaim

\noindent Orbital stability of kinks for $n=1$ is proved in
\cite{HPW}, see also the remark in p.188 \cite{GSS}.
 A  result similar to Theorem 1.1 can be proved for
 traveling kinks.
For $n=2$ and especially for $n=1$
a similar result is open while for $n \ge 4$,
in particular for $n\ge 5$, it should be easier to prove.
When we replace $u_{tt} $ with $iu_t$, we obtain an integrable
Schr\"odinger equation but
 (1.1)  cannot
be treated with the Inverse Scattering Transform, \cite{AS} p. 38.
In the case of the heat equation stronger results than ours are
well known: for $n=1$ see \cite{H};  for $n>1$ see, for $n \ge 4$
\cite{X}, for $n =2,3$ \cite{LX} and for all $n>1$ \cite{Ka}. For
the heat equation the fact that most of the spectrum is strictly
negative is very helpful, while for (1.1) one can think of the
whole spectrum in the imaginary axis.
 In fact for the heat equation when $n=1$ one exploits that 0 is
an isolated eigenvalue and all the rest is strictly negative. For
$n\ge 2$ the point  0 is not isolated
 in the spectrum and so to some extent the
transversal variables complicate the spectral picture. At least in
part this extra difficulty must be purely formal since, for
perturbations localized also in the transversal variables, which
only contributes to their smallness, relaxation to a kink must only
be more likely. This is the view we take in this paper. Since we do
not know how to solve the $n=1$ case we add some extra variables and
exploit the dispersion they provide. The case $n=3$ leads to
equations with a long range nonlinearity.

The equation for the perturbation can be written
$$w_{tt} - \Delta w -w
+3\tanh^2(2^{-\frac 12}x)w+3\tanh (2^{-\frac 12}x)w^2
+w^3=0.\tag 1.4$$
We can rewrite
$$H=- \frac{d^2}{dx^2}- 3 \cosh^{-2} (2^{-\frac 12}x),$$
$$w_{tt} + Hw-\Delta_yw +2w+3\tanh (2^{-\frac 12}x)w^2
+w^3=0.$$ Recall that the eigenvalues of the operator $\dsize
-\frac{\hbar ^2}{\mu}\, \frac{d^2}{dx^2} -\frac{V_0}{\cosh^2
(2^{-\frac 12}x/a)} $ are by  \cite{GK} given by formula
$$\lambda_n = -\frac{\hbar ^2}{\mu a^2}
\left [ \frac 12 \sqrt{\frac{8\mu V_0 a^2}{\hbar ^2}+1} -(n+\frac
12) \right ]^2$$ with $n$ varying among non negative integers so
that $\frac 12 \sqrt{\frac{8\mu V_0 a^2}{\hbar ^2}+1} -(n+\frac
12)>0.$ In our case $H$ has exactly two eigenvalues, given for
$n=0,1$ by
$$\lambda_n = -\frac{1}{2}
\left [ \frac 12
\sqrt{24+1} -(n+\frac 12)
\right ]^2 = -\frac{1}{2} (2-n)^2.$$

In particular $\lambda_0=-2$ and $\lambda_1=-1/2$. Notice
$\lambda_2=0$ is a resonance but not an eigenvalue, that is we have
equality (1.7) with a function bounded but not in $L^2(\Bbb R)$.
That 0 is a resonance is used here crucially, see \S 2 and the proof
of Lemma 12.5. We have
$$\align & (H+2) \cosh^{-2} (2^{-\frac 12}x)=0 \tag 1.5
\\
&\left(H+\frac 12\right) \cosh^{-2} (2^{-\frac 12}x)
\sinh (2^{-\frac 12}x)=0 \tag 1.6
\\
&H \left ( \cosh^{-2} (2^{-\frac 12}x)- 2 \tanh^2 (2^{-\frac 12}x)
\right) =0.
\tag 1.7
\endalign$$
Following Xin \cite{X,LX} and Kapitula \cite{Ka} we write
the solution in the form
$$u(t,x,y)= \th (x-\sigma (t,y))+v(t,x-\sigma (t,y),y) \,, \quad
\int v(t,x,y) \th ^\prime (x)\, dx =0.\tag 1.8$$
In turn, if we set
$$\tilde \phi (x)=\sinh (2^{-\frac 12}x)\cosh ^{-2} (2^{-\frac 12}x)\,
, \, \phi (x)=\frac { \tilde \phi (x) }
{ \|  \tilde \phi ( \cdot)\|_2 } ,\tag 1.9$$
we get
$$v(t,x,y)
=\phi (x) a(t,y) + \psi(t,x,y),\tag 1.10$$ with $ \psi$
corresponding to the continuous spectrum of $H$. In this paper we
prove that locally in space $\psi (t,x,y )=O (t^{-\frac 32+\delta  }
)$ for   $\delta   >0$ small preassigned, $a (t,y )=O (t^{-1} )$ and
$\sigma  (t,y )=O (t^{-\frac 12+\delta  } )$. Probably these
estimates are valid uniformly in space. In the particular case when
for the data in (1.3) we have $w_0(-x,y)=-w_0(x,y)$ and
$w_1(-x,y)=-w_1(x,y)$, where $\sigma  (t,y )=0$ for all $t,y,$ we
are
 able to prove uniform estimates, but we do not include the argument here.

When we plug in (1.1) Ansatz (1.8)
we obtain a system of two NLKG equations, one for $ \psi$
and one for $a$, and a nonlinear wave equation with zero mass  for $ \sigma$.
The nonlinearity
consists  of  pure powers, including quadratic,
 and  null forms.
The equation for $ \psi$ has inhomogeneities that depend on
the longitudinal variable $x$.
If the space dimension $n$ is sufficiently large, then
the theory in Shatah \cite{Sh1}  and dispersion
theory for the linearized operator,
 Weder \cite{We}, lead to the expected result by means of $L^p L^q$
estimates.
If $n$ is small
$L^p L^q$ estimates
are not sufficient to close the inequalities.
The literature
offers as additional
 tools Klainerman's  \cite{K} vector fields  and
Shatah's  \cite{Sh2} normal forms.
 In low dimension $n\le 2$ the two tools are used in conjunction, see
 \cite{D, DFX, GP, Ko, OTT} and
therein  for additional references. The methods  \cite{K,Sh2}
are designed for translation invariant equations.
In \cite{K} a translation invariant NLKG
  is viewed essentially as an ODE in the
radial variable in spacetime. In  \cite{D,DFX} the approach in
\cite{K} is adapted directly to the nonlinear problem.

Since the inhomogeneities in our system depend only on the longitudinal
variable $x$,
we implement the method of  \cite{K,DFX}
(change of  coordinates, normalization of the
 unknowns and  energy estimates on hyperboloids  using Klainerman's
 vector fields) only in the variables $(t,y)$,
while in the $x$ variable the differentiation needed for energy
estimates is by means of the Schr\"odinger operator $H$. The
variables  $(t,y)$ are changed into new ones denoted $T,Y$ while $x$
is left alone. Following   \cite{D,DFX} we consider energy estimates
and enter the information  in the nonlinear system.
 In \cite{DFX} this leads to a
simple system of  ODE's
for the radial variable $T$ plus small integrable perturbations. Then
$L^\infty $ estimates are obtained directly
from the ODE's using standard ODE methods,
for instance standard method of normal forms. Similarly, here
 we obtain a system of one Klein Gordon
equation in $\psi$ with time $T$ and space variable $x$ and  an ODE
in $T$ for $a$, while it is more convenient to think of the equation
for $\sigma $  as a NLW. By variation of parameters and using work
by Weder \cite{We}, we obtain nice estimates for $\psi$,  so that
$\psi$ can be eliminated from the system. Now we have reduced to a
NLW for $\sigma$ and an ODE for $a$. Thanks to  the estimates on
$\psi$, by the Morawetz vector field we get nice estimates for
$\sigma$. Eventually , up to lower order terms, we have in effect
just a closed nonlinear scalar ODE in $a$, with quadratic
nonlinearity in $a$, up to a lower order error. We estimate $a$
applying normal forms as in \cite{DFX}.

We want to point to two features of the problem. The first, which is
certainly fundamental, is the   fact that, both here and in
\cite{Ka},  the problem can be solved only because the
nonlinearities are of a specific form, that is   pure power terms
and, see \cite{Ho, So}, Klainerman's null forms. The second feature,
maybe not as fundamental but important  for our argument,  is that
the endpoints result for $p=1,\infty$
  in Theorem 2.2, more precisely the  dispersive $L^1\to L^\infty $ estimate for linear Klein Gordon equations
  in Corollary 2.3,  are  crucial.
  Now, the  $L^1\to L^\infty $ estimate   to our knowledge is known to hold only with Schr\"odinger
   operators which, like our linearization $H$ defined under (1.4),
   have 0 as a resonance and have
  transition coefficient satisfying $T(0)=1$, see below Lemma 2.1.
  The   $L^1\to L^\infty $ estimate seems essential    because, by the
  dimension 2 of the $Y$ space, at some point we get an undesired $\log T$ growth
  factor in the upper bound on $\psi$, see Lemma 11.2. The estimates
  for $\sigma$ are tight and this $\log T$ term risks to derail all
  the estimates. Fortunately we are able to prove that the crucial term involved in the
  estimate for $\sigma $, see (12.1), does not have the
  additional factor. The proof makes use of the $L^1\to L^\infty $ estimate in Corollary
  2.3. A limitation in the proof is that the  $\delta   >0$ of say $\psi (t,x,y )=O (t^{-\frac 32+\delta  }
)$  affects the size of the neighborhood of $\text{th}$ for initial
data, and not only the constant in the big $O$.

This paper is structured as follows.
In \S 2
we state a framework from Kapitula \cite{Ka}
which is necessary to introduce the modulation equations.
Then we state results in \cite{We} which give
us estimates on groups associated to $H$.
In \S 3 we ``modulate'' following  \cite{Ka}. This leads us to
 a system whose  nonlinearity
is of the right type, in particular displays null form dependence in
$\sigma$. All sections from 4 on are heavily dependent on
\cite{DFX}. In particular we often state without proof formulas and
lemmas which are proved in \cite{DFX}. In \S 5 we introduce
Klainerman's vectorfields and introduce new coordinates. In a
subsection we derive, using Morawetz vectorfield, various formulas
for $\sigma $. In \S 6 we describe the  basic continuity argument
used to prove the main theorem, the rest of the paper consisting in
the proof that inequalities (6.1) imply the improved inequality
(6.2). The crux of the paper starts from \S 7, where we first
restate the system, formulas (7.1-2) and (7.5), and we start a long
list of multilinear estimates. We advise the reader to skip these
estimates at a first reading. In \S 8 we prove a high energy
estimate, Lemma 8.1, using the material in \S 7. At a first reading
we advise the reader to read the statement of Lemma 8.1 and skip the
rest of the section. In \S 9 we rewrite the system emphasizing the
variables $x$ and $T$. In \S 10 we derive improved low energy
estimates for $\Cal A$. In \S 11 we derive dispersive estimates for
$\varPsi$. At the end of \S 11 we exploit the estimates in Lemma
8.6, which are derived from the finite speed of propagation. In \S
12 we prove the estimates for $\varSigma$, stated in Lemmas 12.2 and
12.4. At a first reading we advise the reader to read the statements
of these lemmas, skip all the rest of \S 12 and read the closure of
the estimates in \S 13.

\bigskip
\bigskip
\head \S 2 Spectral decomposition and longitudinal
dispersion \endhead
\bigskip

We start with \S 2 \cite{Ka}.
We  denote by $\langle \rangle _x$ the $L^2$ inner
product in $x$. $\Cal S(\Bbb R)$  is the space of rapidly decreasing and smooth
functions defined in $\Bbb R$.
Then, proceeding as in Lemma 2.1  \cite{Ka}, we
have:
\proclaim{Lemma 2.0} Let $\varphi \in  \Cal S(\Bbb R)$.
For any $U(x,y)$,
$x\in\Bbb R$, $y\in\Bbb R^2$,
let $\langle \varphi , U \rangle _x $
be the function in $y$ obtained by taking
inner product in the variable $x\in \Bbb R$.
Then we have for all integers $k$ and for
all $p \in [1,\infty ]$:
$$\aligned &\| \langle \varphi , U
\rangle _x \|_{W^{k,p}(\Bbb R^2)} \le
\| \varphi \|_{L^{\frac p{p-1}}(\Bbb R)} \| U \|_{W^{k,p}(\Bbb R^3)}
\\
&\| \varphi(x)\langle \varphi, U\rangle _x\|_{W^{k,p}(\Bbb R^3)}
\le
\| \varphi \|_{W^{k,p}(\Bbb R)}
\| \varphi \|_{L^{\frac p{p-1}}(\Bbb R)}
\| U \|_{W^{k,p}(\Bbb R^3)}.\qed
\endaligned$$
\endproclaim

In correspondence
to the spectral decomposition
of $H$, the identity operator in $L^2(\Bbb R)$, in fact
$L^p(\Bbb R)$ for any $p\in [1,\infty ]$,
can be decomposed into
$$\Cal I=P_0+P_1+P_c,\tag 2.0$$
with $P_0$ the projection associated to (1.5), $P_1$ the projection
associated to (1.6) and$P_c$     the projection on the continuous
spectrum of $H$. Since $ \th ^\prime (x) $ and $\phi (x) $ are in
$\Cal S(\Bbb R)$, by Lemma 2.0 projections $P_0$, $P_1$, $P_c$, and
the corresponding splitting (2.0) extend to $W^{k,p}(\Bbb R^3)$ for
any $p\in [1,\infty ]$ and any $k$.

\noindent We collect a few technical facts on operator $H$ needed
later in the proof.  We have: \proclaim{Lemma 2.1} For the
transmission coefficient of $H$, we have $T(0)=1$.
\endproclaim
\noindent By   (2.12) in \cite{We} we have $T(0)= \frac{2a}{2+a^2}$
where $a=\lim _{x\to -\infty } f_1(x,0)$, with $f_1(x,0)$ the Jost
function satisfying $\lim _{x\to +\infty } f_1(x,0)=1.$ If there is
a bounded solution $u$ of $Hu=0$, by elementary ODE arguments $u$
must be a multiple of $f_1(x,0)$. Then the function in (1.7)
appropriately normalized gives $f_1(x,0)$, and  since it is even, we
have $a=1$.\qed
\bigskip
\noindent Lemma 2.1 is very useful because the main result in
 Weder \cite{We}
implies:

\proclaim{Theorem 2.2} Consider for
$\varphi\in C^\infty_0(\Bbb R)$ the operator
$W\varphi = \lim_{t\to +\infty} e^{iHt} e^{i\frac{d^2}{dx^2}t}\varphi.$
Then $W$ extends into a bounded one to one map  from
$W^{k,p}(\Bbb R)$ into itself
$\forall$ $p\in [1,\infty ]$ and $\forall$ $k\in \Bbb Z $,
with image given by $P_c(H)W^{k,p}(\Bbb R)$. \qed
\endproclaim
If 0 was not a resonance, so that $T(0) =0$, or if $T(0) \not =1$,
we would miss $p=1,\infty $ in Theorem 2.2, see \cite{We}. These
endpoints cases are used in Lemma 12.5.

For $P_c(H)$ as in Theorem 2.2 set $B=\sqrt{H+2}P_c(H)$. Since
$Wf(H)P_c(H)= f(-\frac{d^2}{dx^2})W$ for any measurable bounded
function $f$, by Theorem 2.2 $B$ has the following dispersive
properties:

\proclaim{Corollary 2.3}
We have that
$\forall $ $p\in [2,\infty ]$
the operators $\sin ( tB)$ and $\cos ( tB)$
send $P_c(H)W^{k+ 2-\frac 4p,\frac p{p-1}} (\Bbb R) \to
P_c(H)W^{k, p}(\Bbb R)$
with norm bounded by $C_pt^{-\frac 12 + \frac 1p}$. \qed
\endproclaim

\bigskip We have the following fact, see p. 296 \cite{T}:
\proclaim{Lemma 2.4}  $B$ is an elliptic pseudo differential
operator (pdo) of order 1. \qed
\endproclaim
\bigskip

\noindent Finally we have:
\proclaim{Lemma  2.5} Set $B^0=P_c(H)$. If $m$ is even (resp.odd)
 $\forall$
$p\in [2,\infty ]$   (resp.  $p\in [2,\infty [$)  ,
$$\| [B^m, \partial _x ] \colon W^{m,p }
\to L^p \cap L^{\frac p {p-1}} \| < \infty .$$
\endproclaim
For $m=2k$ this follows  from
$$[ (H+2)^k P_c(H),  \partial _x ] = (H+2)^k [ P_c(H),  \partial _x ]
 + [ (H+2)^k ,  \partial _x ] P_c(H).$$
For $m =2k-1$,  $[B^m , \partial _x ]=B [B^{m-1} , \partial _x ]
+[B^{m-1} , \partial _x ]B$ with $[B^{m-1} , \partial _x ]$ bounded
from $  W^{m,p } \to W^{1,p } \cap W^{1,\frac p {p-1}}$ and from $
W^{m-1,p } \to L^p \cap L^{\frac p {p-1}}$. Since $B$ is bounded
from $ W^{m,p } \to W^{m-1,p }$ and from  $W^{1,p } \cap W^{1,\frac
p {p-1}} \to  L^p \cap L^{\frac p {p-1}}$, we are done. \qed

\head \S 3 Modulation \endhead
\bigskip

We still follow Kapitula \cite{Ka}. We need to justify formula
(1.8). For our purposes the following result is sufficient:
\proclaim{Lemma 3.1} Consider $k\in \Bbb N$ and $p\in [1,\infty ]$
such that $ W^{k,p} (\Bbb R^3) \subset C^0 (\Bbb R^3) \cap
L^\infty (\Bbb R^3).$ Then there are positive constants $C$ and
$\epsilon_0$ such that, for any $\epsilon \in ]0, \epsilon_0[$ and
for any $w(x,y) \in W^{k,p} (\Bbb R^3)$ of norm less than
$\epsilon $, there exists a unique pair $(v(x,y),\sigma (y))$
with: $v(x,y) \in W^{k,p} (\Bbb R^3)$; $P_0v=0$; $\sigma (y) \in
W^{k,p} (\Bbb R^2)$;
$$\th (x)+w(x,y)= \th (x-\sigma (y))
+v(x-\sigma (y),y);\tag 3.1$$
 the norm of
$(v(x,y),\sigma (y))$ is less than $C\epsilon.$
\endproclaim
\noindent The proof, that we sketch now, follows
from \cite{Ka} Lemma 2.2.
Assuming (3.1)
we can write
$$w(x,y)=v(x
-\sigma (y),y)- \sigma (y)\int_0^1 \th ^\prime (x
-s\sigma (y)) ds.
\tag 3.2$$
Define
$$F(\sigma, w)(y)= \langle w(\cdot,y), \th ^\prime (\cdot
-\sigma (y)) \rangle _x
+\sigma (y)\int_0^1 \langle \th ^\prime (\cdot
-s\sigma (y)),
\th ^\prime (\cdot
-\sigma (y)) \rangle _x ds.$$
Then
$$F\in C^1( W^{k,p} (\Bbb R^2)\times W^{k,p} (\Bbb R^3)
, W^{k,p} (\Bbb R^2) ),$$
with $F(0,0)=0$ and $F_\sigma (0,0)= \| \th ^\prime \|_2^2\Cal I.$
We have that $P_0v=0$ is equivalent to
 $F(\sigma, w)=0$.
By the Implicit Function Theorem there is a unique continuous map
$\sigma =\sigma (w)$ such that $\sigma (0)=0$ and $F(\sigma
(w),w)=0.$ We then plug this $\sigma $ in (3.2) and we obtain the
desired $v$. \qed

\bigskip
Let us write
$Q_0(f,g)= f_tg_t- \nabla _y f \cdot  \nabla _y g.$
We plug ansatz (1.8) into equation (1.1) and
get, renaming $x-\sigma$ by $x$,
$$\aligned & Q_0(\sigma ,\sigma ) \th ^{\prime \prime } (x) -
  \th ^\prime (x)\square \sigma
- v_{x}(x)\square \sigma +(3\th  v^2
+v^3)(x)
\\
&+(\square v +2v+3( \th ^2 -1)v) (x) + Q_0(\sigma ,\sigma )v_{xx}
-2Q_0(\sigma ,v_{x}  )=0.
\endaligned \tag 3.3$$
In the frame associated to (2.0), if we write
$v(t,x,y)=P_1v(t,x,y)+P_c v(t,x,y)
=\phi (x) a(t,y) + \psi(t,x,y)$, equation (3.3) splits into
a system formed by (3.4-6) below. We will set $\square _y=
\partial _{tt} -\Delta _y.$ We consider first $P_c(3.3)$
to obtain:

$$\aligned
&(\square _y+B^2)\psi-P_c( \psi _{x}+ a \, \phi ^{\prime })\square _y \sigma
+Q_0(\sigma ,\sigma )P_c\psi _{xx}- 2 Q_0(\sigma ,P_c\psi _{x}  )= \\&
 = -P_c (F_2+F_3) -Q_0(\sigma ,\sigma )P_c  \th ^{\prime \prime }
+2 Q_0(\sigma ,a) P_c\phi   ^{\prime }
-a
Q_0(\sigma ,\sigma )P_c  \phi   ^{\prime \prime}
\\
&F_2(\psi,a)= 3\th \, (a\phi +\psi )^2
\\
& F_3(\psi,a)
=(a \, \phi +\psi )^3.\endaligned \tag 3.4$$
We consider $\langle (3.3) , \phi \rangle _x$  (from
now on we omit the subscript $x$
and write simply $\langle \rangle $) and obtain
$$\aligned
&(\square _y+\frac 32)a =
-G_2-G_3-\langle \psi ,\phi  ^{\prime \prime} \rangle \square \sigma
-Q_0(\sigma ,\sigma ) \langle \th ^{\prime \prime }  , \phi   \rangle
\\&
-2 Q_0(\sigma , \langle \psi  , \phi ^{ \prime} \rangle  )  -
Q_0(\sigma , \sigma ) \langle \psi  , \phi  ^{\prime \prime} \rangle   +a
Q_0(\sigma ,\sigma )  \| \phi    ^{\prime} \| ^2_2
\endaligned \tag 3.5$$ where $G_j=\langle F_j ,
\phi \rangle .$
Notice that
$2Q_0(\sigma ,a ) \langle \phi  , \phi ^\prime  \rangle =0$.
We consider $\langle (3.3) ,  \th ^{\prime \ } \rangle $ and obtain
$$\aligned  &\square _y\sigma \,  ( \|  \th ^{ \prime }  \| ^2_2
- \langle \psi  ,   \th ^{\prime \prime }  \rangle
+ a \langle   \phi    ^{ \prime } ,   \th ^{\prime }  \rangle )=H_2+H_3+\\& +
2 Q_0(\sigma , \langle \psi  ,  \th ^{\prime \prime }  \rangle  )
-2Q_0(\sigma , a )\langle \phi  ^{ \prime } ,  \th ^{ \prime }  \rangle  +
Q_0(\sigma , \sigma )\langle \psi  ,
 \th ^{\prime \prime  \prime }  \rangle  \endaligned
\tag 3.6 $$
where $H_j=\langle F_j , \th
    ^\prime   \rangle .$ In particular
$$ H_2(\psi,a)
=6a \langle  \psi,
 \th \,  \phi  \,  \th ^{ \prime } \rangle
+3\langle  \th   \,   \psi^2,  \th ^{ \prime }  \rangle$$ where
$3a^2 \langle  \th  \,   \phi^2,  \th ^{ \prime }  \rangle =0$
because $\phi^2 \th ^{ \prime } $ is even while $ \th  $ is odd.
Notice that the nonlinearities contain $\sigma$ only inside   null
forms. We will prove:

\proclaim{Theorem 3.2} Under the hypothesis of Theorem 1.1 we have
that for any fixed $\delta >0$ there is an $\epsilon _0 =\epsilon
_0(\delta )>0$ such that for initial data as in Theorem 1.1 we have
for $T(t,y)=\sqrt{1+|t^2-|y |^2| }$

$$| \psi (t,\vec{x} )| T^{\frac 32-\delta }(t,y)+ |a (t,y )| T (t,y)+ |\sigma  (t,y )| T ^{\frac 12-\delta }(t,y)
 < \infty .$$
 \endproclaim

\noindent Having completed the set up, in the rest of the paper we
borrow heavily from \cite{DFX}.

\bigskip
\bigskip
\head \S 4 Short term behaviour \endhead
\bigskip

We first look at equation (1.4) with initial conditions (1.2).
In analogy to Proposition 1.1.4 \cite{DFX} and
Proposition 1.4 \cite{D}, we have:
\proclaim{Proposition 4.1} Let $T_0>2 K$, $s \in \Bbb N$,
$s \ge 3$.
There are $\epsilon_0>0$ and $C>0$
such that for any $\epsilon$ with $0<\epsilon <\epsilon_0$
and for any
$(w_0,w_1)$ in a ball of radius $\epsilon$ in $H^{s+1} (\Bbb R^3)
\times H^{s} (\Bbb R^3)$ and supported in
$\{ (x,y) \colon
|x|^2+ |y|^2\le K^2\} $, problem (1.4) (1.2) has a unique solution
defined in
$$\{(t,y): t\ge 0,\,(t+2K)^2-|y|^2\le T_0^2 \}, \tag 4.1$$
continuous in time with values in $H^{s+1} (\Bbb R^3)$
and $C^1$ with values in $ H^{s} (\Bbb R^3)$.
Moreover
$$\sum_{|\alpha |\le s} \int_{\Bbb R^3}
\big | (\partial^\alpha w)(\sqrt{T_0^2 +|y|^2},x, y)
\big |^2dx \, dy \le C^2 \epsilon^2, \tag 4.2$$
and the restriction of $w(t,x,y)$ and its derivatives
on the hypersurface
$(t+2K)^2-|y|^2= T_0^2$ is supported inside
$\{ \sqrt{|x|^2+ |y|^2}\le \frac{T_0^2-K^2}{2K} \}.$
\endproclaim
We sketch the proof. By taking $\epsilon_0$ small, there is a
unique solution $w(t)\in L^\infty ( [0,t_0),H^{s+1} (\Bbb R^3))
\cap C^{0,1} ( [0,t_0),H^{s} (\Bbb R^3)) $ for any preassigned
$t_0$, Theorem 6.4.11 \cite{Ho}, which is also in
 $C^{0}( [0,t_0),H^{s+1} (\Bbb R^3)) \cap C^1( [0,t_0),H^{s} (\Bbb R^3))
$,   \cite{So} p.18.
For $t_0=\frac{T_0^2-3K^2}{2K}$,
the intersection
of (4.1) with the support of $w(t,x, y)$ is contained
in $\{ t\le t_0\} \cap \{ (|x|^2+ |y|^2)^{\frac 12}\le t_0+K\} .$
(4.2) can be obtained
by a trace theorem from the standard energy inequality.
\bigskip
\head \S 5 Change of coordinates
and
Klainerman's vector fields
\endhead
\bigskip

Having disposed of the solution in the region below (4.1),
we consider system (3.4-6) in the region above (4.1).
Following \cite{K,D,DFX},
we replace coordinates
$(t,x,y)\to (T,x,Y)$, we introduce in this new
system of coordinates Klainerman's
vector fields, only for variables
$(t,y)$, and we rewrite
(3.4-6) in the new coordinate system.

\noindent Set $t+2K=T\cosh |Y|$, $y=TY \frac{\sinh |Y|}{|Y|}$.
Introduce Klainerman vector fields
$$\aligned Z_0&=y_1\partial_{y_2}-y_2\partial_{y_1}
=Y_1\partial_{Y_2}-Y_2\partial_{Y_1}
\\
Z_1&= (t+2K) \partial_{y_1}+y_1\partial_{t}
\\
&=\left[\frac{Y_1^2} {|Y|^2}
+\frac{Y_2^2\cosh |Y|} {|Y|\sinh |Y|}\right]
\partial_{Y_1}
+\frac{Y_1Y_2} {|Y|^2}
\left[1-\frac{|Y|\cosh |Y|} {\sinh |Y|}
\right]\partial_{Y_2}
\\
Z_2&= (t+2K) \partial_{y_2}+y_2\partial_{t}
\\
&= \frac{Y_1Y_2} {|Y|^2}
\left[1-\frac{|Y|\cosh |Y|} {\sinh |Y|}\right]
\partial_{Y_1}
+\left[\frac{Y_2^2} {|Y|^2}
+\frac{Y_1^2\cosh |Y|} {|Y|\sinh |Y|}\right]
\partial_{Y_2}\\
Z_3&=\partial _ T
.
\endaligned$$

We denote by $\Cal E^r$ the set of functions $(T,x,Y) \to c(T,x,Y)$
such that for any multi-index $I$ and integer $k$ there is a
constant $C_{I,k}$ such that
$$| \partial_x^k  Z^I c(T,x,Y) |\le C_{I,k}T^{-m}e^{r|Y|}$$
where $m$ is the number of $Z_3$ factors inside $Z^I$.

By Lemma 1.2.2 \cite{DFX}  we have, with $b ^{(\ell )}_j(Y) \in \Cal
E^0,$ $$\partial_{Y_j}= b ^{(1)}_j(Y) Z_1+ b^{(2)}_j(Y) Z_2 .\tag
5.1$$ Fix spherical coordinates $ Y=R(\cos \theta,\sin\theta)$. We
consider
$$\partial_R= \frac Y{|Y|} \cdot \nabla_Y \,, \quad
\partial_\theta =-Y_2 \partial_{Y_1}+Y_1 \partial_{Y_2}=\tanh (R)
\big [- \frac {Y_2}R Z_1+\frac {Y_1}R Z_2 \big ] .$$ Set
$$\aligned
&\partial_t^2- \Delta_y =P+\frac 2T \partial_T\, ,
\quad P= \partial_T^2-\frac 1{T^2} \Delta_{hyp} \quad \text{where}\\
&  \Delta_{hyp} =\partial_R^2+\frac{\cosh R} {\sinh R}
\partial_R+\frac 1 {\sinh^2R} \partial_\theta^2 =
Z_1^2+Z_2^2-Z_0^2 .
\endaligned \tag 5.2$$
Next , we consider
 for $s\in \Bbb N$
the   spaces $H^s_{ Y}$  with norm, for $I =  ( I_0, I_1, I _2)$ and
$|I |=|I _{0}|+|I _{1}|+|I _{2}|, $

$$\| f  \| _{H^s_{ Y}}^{2}=\sum _{|I |\le s} \int _{\Bbb R^2} |
Z^I f(Y)|^2\sinh (R) dRd\theta . \tag 5.3$$ We set $L^2_Y= H^0_Y$.
We consider
 for  $s \in \Bbb N$ the
Sobolev spaces $H^s_{xY}=H^s_{ Yx}$ with norm

$$\| f  \| _{H^s_{ xY}}^{2}=\sum _{|I |+m\le s} \int _{\Bbb R^2} |
\partial _x^mZ^I f(Y)|^2\sinh (R) dRd\theta dx. \tag 5.4$$

\noindent We denote $H^0_{xY}$ by $L^2_{xY}$ or $L^2.$

\noindent

We now consider the null form $Q_0$ in the new coordinates. By
elementary computation:
$$Q_0(f,g)
=f_Tg_T-\frac 1{T^2}f_Rg_R
-\frac 1{T^2}\frac{1} {\sinh^2 (R)} f_\theta g_\theta .$$

\noindent As in \cite{K} we set
$( \psi, a,\sigma)
=\frac{1}{T}( \varPsi,\Cal{A},\varSigma).$
Notice that the above normalization leads to weaker results
than in \cite{DFX}.
Set
$$\aligned & \Cal Q _0(f ,g )
=( f_{T}-\frac 1T f) ( g_{T}-\frac 1T g) - \frac{  f_{R}g_{R}
}{T^{2}} -\frac{f_{\theta}g_{\theta}} {T^2\sinh^2 (R)} .\endaligned
\tag 5.5$$ With the above changes of the independent variables and
of the unknowns, system (3.4-6) becomes

$$\align
(P+B^2)\varPsi -&\frac 1T P_c(  \varPsi_{x}+ \Cal A
 \, \phi _x)P  \varSigma
+\frac 1{T^2}\Cal Q_0(\varSigma , \varSigma )P_c\varPsi _{xx} -
 \frac 2T\Cal Q_0( \varSigma ,P_c\varPsi _{x}  )=\tag 5.6 \\
 = &- \frac {P_c} TF_2(\varPsi , \Cal A \phi (x))
-\frac 1T
Q_0( \varSigma ,\varSigma )P_c\th ^{\prime\prime }
+\frac 2T \Cal Q_0(\varSigma ,\Cal A) P_c\phi ^{\prime }  -\\&
-\frac  {P_c}{T^2}F_3(\varPsi , \Cal A \phi (x))
   -\frac 1 {T^2}\Cal A\Cal
Q_0( \varSigma ,\varSigma )P_c  \phi     ^{\prime \prime }
\\
(P+\frac 32)\Cal A =& -\frac 1TG_2(\varPsi , \Cal A)-\frac
1{T^2}G_3(\varPsi , \Cal A)- \frac 1{T}\Cal Q_0( \varSigma
,\varSigma )\langle \th ^{\prime\prime }  , \phi  \rangle - \frac 1T
\langle \varPsi ,\phi  ^{\prime \prime } \rangle P \varSigma \tag
5.7\\& -\frac 2T \Cal Q_0( \varSigma, \langle  \varPsi  , \phi  ^{
\prime } \rangle  )  - \frac 1{T^2}\Cal Q_0( \varSigma ,\varSigma
)\langle \varPsi  , \phi  ^{\prime \prime } \rangle +\frac
1{T^2}\Cal A \Cal Q_0( \varSigma ,\varSigma ) \| \phi    ^{ \prime }
\| ^2_2
\\
(
 P\varSigma ) \,  ( \| \th ^{\prime } \| ^2_2 &
- \frac 1T\langle \varPsi  ,  \th ^{\prime\prime } \rangle + \frac
1T \Cal A \langle   \phi    ^{ \prime } ,  \th ^{\prime } \rangle )
=\frac 1TH_2(\varPsi , \Cal A) +\frac 2T \Cal Q_0( \varSigma ,
\langle \varPsi  , \th ^{\prime\prime } \rangle) - \tag 5.8\\&
-\frac 2T \Cal Q_0( \varSigma , \Cal A)
 \langle \phi  ^{ \prime }  , \th ^{\prime } \rangle
+\frac 1{T^2}H_3(\varPsi , \Cal A) +
 \frac 1{T^2}\Cal
Q_0( \varSigma ,\varSigma )
\langle \varPsi  ,  \th ^{\prime\prime \prime  } \rangle  .\endalign
 $$

\noindent By Proposition 4.1 $w$ is smooth in (4.1) with  bounds on
the first few derivatives. The same statement holds for
$(\varPsi,\Cal{A},\varSigma)$ by Lemma 2.0. We will prove:
\proclaim{Theorem 5.1} Let $s_0 =23$. Fix  a small positive number
$\delta  >0$. Let $\epsilon_0 >0$ and let $\epsilon $ be any number
$\epsilon \in ]0,\epsilon_0 [ $. Consider smooth initial data for
(5.6-8)
$$\aligned & \varPsi (T_0)= \varPsi_{(0)}
\,,\, \partial_T\varPsi (T_0)= \varPsi_{(1)}\,,\,
\Cal{A} (T_0)=\Cal{A}_{(0)}
\,,\, \partial_T\Cal{A} (T_0)=\Cal{A}_{(1)},
\\
&\varSigma (T_0)= \varSigma_{(0)}
\,,\, \partial_T\varSigma (T_0)= \varSigma_{(1)}
\endaligned
\tag 1$$ which are the traces on $T=T_0$ of the solutions with
smooth initial data provided by Proposition 4.1. Suppose we have
estimates
$$\| (\varPsi_{(0)},\Cal{A}_{(0)},\varSigma_{(0)})\|_{H^{s_0}}
+\| (\varPsi_{(1)},\Cal{A}_{(1)},\varSigma_{(1)})\|_{H^{s_0-1}} \le
\epsilon.\tag 2$$  Then it is possible to choose $\epsilon_0$ such
that system (5.6-8) with initial conditions (1) satisfying (2) has a
unique solution in $C^{0}([T_0, \infty), H^{s_0})\cap C^{1}([T_0,
\infty), H^{s_0-1})$ and such that we have $$(T^{\frac 12-\delta  }
\varPsi (T)  , \Cal{A} (T) , T^{-\frac 12-\delta  } \varSigma  (T) )
\in L^\infty _{xY}.$$
\endproclaim
The hypothesis that the initial data are traces of solutions from
Proposition 4.1 is used later when we need to exploit the finite
speed of propagation of  the solution $w$ of (1.4).
\bigskip
\noindent In the region of existence the solution $
(\varPsi,\Cal{A},\varSigma)$
is smooth.
For $q\in \Bbb N$, set
$$\varPsi _q=( Z^IB^m \varPsi )_{m+|I|\le q}\, , \,
\Cal A _q=( Z^I \Cal A  )_{|I|\le q} \, , \, (\varPsi
,\Cal{A})_q=(\varPsi _q,\Cal{A}_q).\tag 5.9$$ Similarly set  $
\varSigma _q=( Z^I  \varSigma  )_{|I|\le q}.$ We state:
\proclaim{Lemma 5.2} We have: \roster

\item $[Z_j,P]=0$ for $j=0,1,2$,see \cite{DFX} Lemma 1.2.3;

\item $[Z_j,\partial _T]=0$ for all $j$, see \cite{DFX} (1.2.10);

\item $[Z_j,\Delta _{hyp}]=0$ for all $j$, from the above two claims

\item $[\partial _T,P]=[\partial _T,- \frac
1{T^2}\Delta _{hyp} ]= \frac 2 {T^{3}}\Delta _{hyp}$.\qed
\endroster
\endproclaim
\bigskip
\noindent Set $[Z,\Cal B ] (u,v)= Z\Cal B (u,v)-\Cal B (Zu,v)-\Cal B
(u,Zv) $ for any vectorfield $ Z$ and  bilinear form $\Cal B$. By
elementary computations, using formulas $Q_0(f,g)=\square
(fg)-g\square f-f \square g$ and $\Cal Q_0(f,g)=P (fg)-gP f-f P g$
we get: \proclaim{Lemma 5.3 } Consider a multiindex $I$. In view of
Lemma 5.2 we can write $Z^I=Z_3^k Z^{I^\prime } $, where $
Z^{I^\prime }$ does not contain $Z_3=\partial _T$ among its factors.
We also set $\Cal Q_1= T^2 \Cal Q_0.$ The we have: \roster

\item
$[Z^{I^\prime },\Cal Q_0 ]=0$;

\item $[ Z^{I^\prime },\Cal Q_1 ]=0$;

\item $[ T\partial _T,\Cal Q_1 ]=0$;

\item $[ T\partial _T,\Cal Q_0 ]=-2\Cal Q_0  $. \qed
\endroster
\endproclaim

\bigskip
\head \S 5.1 The use of Morawetz vectorfield \endhead

 Lemmas 5.5 and 5.6 are crucial for the dispersive estimates
 for $\sigma $ derived in \S 12. In particular we will first derive Lemma 5.5, which is
later used in Lemma 12.2. Then we will complement Lemma 5.5 with
Lemma 5.6, obtained through a simple computation and leading to the
estimates of Lemma 12.4.

To discuss Lemma 5.5 we introduce the Morawetz vectorfield
$$K_0= (t^2+|y|^2) \partial_t+2ty\cdot \nabla_y +t
=T^2 \cosh ( |Y|) \partial_T+ T \sinh ( |Y|)\partial_R+T \cosh ( |Y|).$$
We  write
$$K_0\frac{\varSigma}{T }
=\frac{1}{T } \Cal K \varSigma \quad \text{
with}\quad
\Cal K = T^2 \cosh ( R) \partial_T+T \sinh ( R)
\partial_R.$$
 We have:

\proclaim{Lemma 5.4} We have
$$\aligned
&\sinh (R)\, (\Cal Ku)\, (Pu)
=\partial_T\Cal P^0+\partial_R\Cal P^1+\partial_\theta\Cal P^2
\endaligned$$
 where
$$\aligned  \Cal P^0=&\frac 12\sinh ( R)\cosh  ( R) \big [
T^2u_T^2 + 2T\tanh ( R)u_Tu_R
+ u_R^2 +\frac {u_\theta^2}{\sinh ^2( R)}\big ]
\\
\Cal P^1=&\frac 12 \big [ - T\, u_T^2-\frac 2 {
\tanh  ( R)}\, u_Tu_R
- \frac 1T\,  u_R^2
+\frac { u_\theta^2}{T \sinh ^2 ( R)} \big ] \sinh ^2 ( R) \\ \Cal P^2=&
-T^{-1}\, u_R u_{\theta  }-\tanh ^{-1} ( R)\, u_T u_{\theta  }.
\endaligned$$
\endproclaim
PROOF. With a  change of coordinates $(t,y)\to (T,Y)$ and for
$f=T^{-1}u$, we get
$$dtdy\,  (K_0f) \,  \square _{ty} f=dT \,  dR \,  d\theta \,  \sinh  ( R)
 \, T  (T\cosh ( R) u_T+ \sinh ( R) u_R) Pu.$$
We have
$$ T^2\sinh  ( R)
\cosh  ( R)u_T u_{TT}= (\frac {T^2}2\sinh  ( R)
\cosh ( R) u_T^2 )_T- T\sinh ( R)
\cosh  ( R)u_T^2 .\tag 1$$
We have
$$ \aligned T \,\sinh ^2( R) u_R  u_{TT}&= (T \,\sinh ^2( R) u_R  u_{T})_T-
\sinh ^2( R) u_R  u_{T}-T \,\sinh ^2( R) u_{RT}  u_{T} \\& =
(T \,\sinh ^2 ( R)u_R  u_{T})_T-
\sinh ^2 ( R) u_R  u_{T}-\\& - (\frac T2 \,\sinh ^2 ( R)u_{T}^2)_R +
 T \,\sinh  ( R)  \,
\cosh R u_T^2
.\endaligned \tag 2$$
We have
$$ \aligned & -\sinh  ( R)
\cosh ( R)u_T\,  u_{RR}=  -(\sinh ( R) \cosh ( R)\, u_T
u_{R})_R+\\&+ \frac 12(\sinh( R) \cosh ( R)\, u_R^2)_T+ (\cosh ^2(
R) +\sinh^2  ( R))  u_T u_{R} .
\endaligned \tag 3$$
We have
$$-\frac 1T \,\sinh ^2( R)\,  u_R u_{RR}=
-\frac 12(\frac 1T \,\sinh ^2( R)\,  u_R ^2)_R+
\frac 1T \,\sinh ( R)\,\cosh  ( R)\,  u_R ^2
 .\tag 4$$
Another term is
$$ \aligned & -
\cosh ^2( R)u_T u_{R} .
\endaligned \tag 5$$
Another term is
$$-\frac 1T\,\sinh ( R)\, \cosh ( R) \,  u_R^2
 .\tag 6$$
We have
$$-\tanh ^{-1} ( R)\, u_T u_{\theta \theta }=
-(\tanh ^{-1} ( R)\, u_T u_{\theta  })_ \theta  +\frac 12
(\tanh ^{-1} ( R)\,  u_{\theta  }^2)_ T  .\tag 7$$
We have
$$-T^{-1}\, u_R u_{\theta \theta }=
-(T^{-1}\, u_R u_{\theta  })_ \theta  +\frac 12
(T^{-1}\,  u_{\theta  }^2)_ R  .\tag 8$$
Adding up  the numbered formulas
and simplifying we obtain Lemma 5.4.\qed
\bigskip
By Lemma 5.4 we have
$$\aligned
\partial_T(\frac{\Cal P^0}{\cosh  ( R)})
+\partial_R(\frac{\Cal P^1}{\cosh  ( R)})
+\partial_\theta (\frac{\Cal P^2}{\cosh  ( R)})+
\frac {\tanh ( R)}{\cosh ( R)} \Cal P^1= \tanh ( R) \,  \Cal K u  \, Pu.
\endaligned$$
By the last formula we have:

\proclaim{Lemma 5.5} Suppose $P u=F(T,Y)$.
Then, for $\Cal P^j $ as above, we have:
$$\aligned &
\partial_T(\Cal P^0/\cosh ( R))+\partial_R(\Cal P^1/\cosh ( R))
+\partial_\theta (\Cal P^2/\cosh ( R))= \\& = \sinh ( R) \, (T^2
\partial_Tu+T\tanh ( R) \, \partial_Ru)F- \frac {\tanh ( R)}{\cosh
( R)} \Cal P^1 .\qed
\endaligned$$
\endproclaim
\bigskip

To get the desired dispersion for $\sigma $  we will need to
supplement Lemma 5.5 with  the following elementary computation:

\proclaim{Lemma 5.6} Suppose $P u=F(T,Y)$. Then we have
$$\aligned   \sinh ( R) \, T^2 u_T Pu&=  \left [ \frac{\sinh ( R)}2
(T^2u_T^2+ u_R^2+\frac {u_{\theta  }^2} {\sinh ^2( R)} )\right ] _T -\\&-
(\sinh ( R) \,  u_R u_{T  })_R- (\sinh ^{-1}( R) \,  u_\theta u_{T  })_\theta
-T\sinh ( R) \, u_T^2.
\endaligned$$
\endproclaim
The proof follows from a direct computation. Here notice
$$\aligned &  \sinh ( R) \, T^2 u_Tu_{TT}=
\frac 12 ( \sinh ( R) \,  T^2u_T^2)_T-T\sinh ( R) \, u_T^2, \\&
- \sinh ( R) \,  u_Tu_{RR}=-( \sinh ( R) \,  u_Tu_R)_R +\cosh ( R) \,  u_Tu_R
+\frac 12 ( \sinh ( R) \,  u_R^2)_T, \\&
-\sinh ^{-1}( R) \,  u_Tu_{\theta \theta }=-
(\sinh ^{-1}( R) \,  u_Tu_{\theta  })_\theta  +\frac 12
( \sinh ^{-1}( R) u_{\theta  } ^2)_T.
\endaligned$$

\bigskip
\bigskip
\head \S 6 Energy and a continuity argument \endhead
\bigskip

\noindent We define, for appropriate $f(T,Y)$,
$$\tilde E(T,f)
=\int_{\Bbb R^2} \left[ | f_{T}|^2
+\left| \frac{\partial_R}{T} f \right|^2
+\big | \frac{1} {T\sinh ( R)}
\partial_\theta f \big |^2
\right] \sinh ( R) dR d\theta .$$ Next define
$$\aligned & E_{(1)}(T,\varPsi)=
\int_{\Bbb R} \tilde E(T,\varPsi) dx
+\int_{\Bbb R^3} |B\varPsi |^2  \sinh ( R)\,  dR d\theta\,dx
\\
&E_{(2)}(T,\Cal{A})=
 \tilde E(T,\Cal{A})
+ \int_{\Bbb R^2} \frac 32|\Cal{A} |^2  \sinh ( R)\, dR d\theta .
\endaligned$$
We call $E(T,(\varPsi,\Cal{A} ))=E_{(1)}(T,\varPsi)+E_{(2)}(T,\Cal{A})$.
We have, see \cite{DFX} Lemma 2.1.1:

\proclaim{Lemma 6.1}
We have
$$\aligned & \frac d{dT}E_{(1)}(T,\varPsi)\le  -
2\Re \int_{\Bbb R^3} \bar \varPsi_{T}
(P \varPsi +B^2 \varPsi)  \sinh ( R)\, dR d\theta\,dx \\
 & \frac d{dT}E_{(2)}(T,\Cal{A})\le
-2\Re \int_{\Bbb R^2} \bar{\Cal{A}}_{T}
(P\Cal{A} +\frac 32\Cal{A})  \sinh ( R)\, dR d\theta
\\ & \frac d{dT}\tilde E(T,\varSigma)\le
-2\Re \int_{\Bbb R^2} \bar \varSigma_{T}
P \varSigma  \sinh ( R)\, dR d\theta .
\endaligned$$
\endproclaim

\noindent Next consider, see (5.5), $E_q(T)=E(T, (\varPsi,\Cal{A})_q)=
E_{(1)}(T,\varPsi _q) +E_{(2)}(T,\Cal{A} _q)$ where
$$E_{(1)}(T,\varPsi _q) =\sum_{m+|I|\le q}E_{(1)}(T,B^mZ^I\varPsi )\,
, \, E_{(2)}(T,\Cal{A} _q)=\sum_{|I|\le q}E_{(2)}(T,Z^I \Cal{A}).$$
Similarly set
$$\tilde E_q
(T)=\tilde E(T,\varSigma _q)= \sum_{|I|\le q}\tilde E(T,Z^I
\varSigma ).$$ We fix integers  $N=22$ and
 $N^\prime =14.$
$N$ and $N^\prime $  are chosen so that $N\ge N^\prime +8$ and $
N^\prime \ge [\frac N2]+3$ and $N^\prime \ge 8$. Next we fix a
small number $\delta >0$. We then fix $p\in (2,\infty )$ so that
$\delta \ge\frac 2p$ (we pick $p\neq \infty $ because pdo's are
not well behaved in $L^\infty $). Next, we suppose that in an
interval $[T_0, T^\ast [$ we have
$$\aligned & E(T_0,(\varPsi,\Cal{A})_N)+ \tilde
E(T_0, \varSigma _N) < \epsilon^2
\\
&\sup_{T\in [T_0, T^\ast [} ( T^{\frac 12-2\delta} \| \varPsi
_{N^\prime }\| _ {L^2_YL^p_x}+ T^{\frac 12-\delta }\tilde  E^{\frac
12}(T,\varSigma_{N^\prime }) ) \le \mu ^\prime \epsilon \\ &
\sup_{T\in [T_0, T^\ast [} (   \| \varPsi _{N^\prime }\| _
{L^2_YL^2_x}+ \| \partial _T\varPsi _{N^\prime }\| _ {L^2_YL^2_x}+
\| \Cal{A} _{N^\prime }\| _ {L^2_Y } +\| \partial _T\Cal{A}
_{N^\prime }\| _ {L^2_Y } ) \le \mu ^\prime \epsilon .
\endaligned \tag 6.1$$
The crux of the proof consists in the following continuity argument.
We need to show that we can choose  $\mu ^\prime $ and $\epsilon
(\mu ^\prime )$ so that
 for any $\epsilon \in ]0, \epsilon (\mu ^\prime )[$,
then (6.1) implies
$$\aligned & \sup_{T\in [T_0, T^\ast [} ( T^{\frac 12-2\delta} \| \varPsi
_{N^\prime }\| _ {L^2_YL^p_x}+ T^{\frac 12-\delta }\tilde  E^{\frac
12}(T,\varSigma_{N^\prime }) ) \le \frac{\mu ^\prime}{2} \epsilon \\
& \sup_{T\in [T_0, T^\ast [} (   \| \varPsi _{N^\prime }\| _
{L^2_YL^2_x}+ \| \partial _T\varPsi _{N^\prime }\| _ {L^2_YL^2_x}+
\| \Cal{A} _{N^\prime }\| _ {L^2_Y } +\| \partial _T\Cal{A}
_{N^\prime }\| _ {L^2_Y } ) \le \frac{\mu ^\prime}{2} \epsilon .
\endaligned \tag 6.2 $$
Once we have that $(6.1) \Rightarrow (6.2)$ we are done. Indeed,
suppose
  $T^\ast <\infty$.
In the region $T\in [T_0, T^\ast [$, by the finite speed of
propagation the support of $w(t,x,y)$ is a bounded set. For any
point $P_0=(t_0,x_0,y_0)$, on the hypersurface $T=T^\ast $ the
behavior of $w(t,x,y)$ depends only on the values of $w(t,x,y)$ in
the cone defined by $t_0-\epsilon _1<t <t_0$, for any $\epsilon _1
>0$ and $|(x,y)- (x_0,y_0)|\le t_0-t.$ For $\epsilon _1 $ small, the
cone is in the region $T\in [T_0, T^\ast [$. On the cone, (6.1)
implies that a large number of derivatives of $w(t,x,y)$ is bounded,
and so
 $w(t,x,y)$ does not blow up at $P_0$. Therefore we can conclude that
there is a $T_1> T^\ast $ such that $w(t,x,y)$
 can be extended in a
$[T_0, T_1 [$. By ($6.2$) we can assume (6.1)
is valid in this larger region. This implies we
 can choose $T^\ast =\infty $.
\bigskip

\head \S 7 Various   inequalities \endhead

We start by rewriting equations (5.3-5). Set
$$\Omega = \|   \th ^{ \prime } \| _2^2
- \frac 1T \langle \varPsi  ,  \th ^{\prime \prime } \rangle + \frac
1T \Cal A \langle   \phi   ^{ \prime }  ,    \th ^{ \prime } \rangle
.\tag 7.1$$

\proclaim{Lemma 7.1} There is a a constant $C(\mu ^\prime )$ such
that $\forall $ $T\in [T_0,T^\ast [$ we have
$$\aligned &
\|  \varSigma _{N^\prime -1 }(T) \| _{L^\infty _Y} \le C(\mu ^\prime
) \epsilon T^{\frac 12+\delta }\, , \, \| \partial _T \varSigma
_{N^\prime -2 }(T) \| _{L^\infty _Y} \le C(\mu ^\prime ) \epsilon
T^{-\frac 12+\delta }\\ & \|  \varPsi _{N^\prime -2 }(T) \|
_{L^\infty _YL^p_x} \le C(\mu ^\prime ) \epsilon T^{-\frac 12+\delta
}\, , \, \| \Cal A _{N^\prime -2 }(T) \| _{L^\infty _Y} \le C(\mu
^\prime ) \epsilon  . \endaligned$$
\endproclaim
\noindent By
  Sobolev Embedding Theorem, see \cite{DFX} Corollary 2.2.4,
and (6.1).\qed

\bigskip
By   Lemma 7.1, $\forall \, T\in [T_0, T^\ast [$ we have $\Omega
\approx \|   \th ^{ \prime } \| ^2_2 $ for $\epsilon \ll 1$. We will
denote schematically $L( \varPsi ,\Cal A )=\lambda
 \Cal A+\mu \langle \varPsi , \psi \rangle $ (or $L(  \varPsi
,\Cal A )=\lambda
 \Cal A \psi (x)+\mu   \varPsi   $)
for $\lambda $ and $\mu $ constants with $|\lambda |+|\mu |\le C$
for a fixed $C$, and for $\psi (x)$ some Schwartz function. For $L$
changing from place to place, we can write schematically

$$\aligned &
 P\varSigma
= \frac{\Cal H }{\Omega   } \, , \quad \Cal H=
\frac{1}{T}H_2(\varPsi , \Cal A) +\frac{1}{  T }\Cal Q_0( \varSigma
, L(  \varPsi ,\Cal A )) +\\& +\frac 1{  T^2}H_3(\varPsi , \Cal A) +
 \frac 1{  T^2}\Cal
Q_0( \varSigma ,\varSigma ) L(  \varPsi ,\Cal A )  .\endaligned \tag
7.1
 $$
Next, let us write schematically

$$\aligned &
(P+\frac 32)\Cal A =  -\frac 1TG_2(\varPsi , \Cal A)-\frac
1{T^2}G_3(\varPsi , \Cal A)- \frac 1{T}\Cal Q_0( \varSigma
,\varSigma )\langle \th '', \phi \rangle   \\& - \frac 1T \langle
\varPsi ,\phi ^{\prime \prime } \rangle \frac{\Cal H }{\Omega   }
-\frac 2T \Cal Q_0( \varSigma, L( \varPsi ,\Cal A ) )  - \frac
1{T^2}\Cal Q_0( \varSigma ,\varSigma )L( \varPsi ,\Cal A )
  .\endaligned \tag 7.2
 $$
We rewrite now the equation for $\varPsi$. We consider a symmetric
matrix with entries
$$\aligned & r_{xx} =\frac 1{T^2} \Cal Q_0( \varSigma , \varSigma ) \quad ,
  \quad
r_{xR}=r_{Rx} =\frac 2{T^3}  \varSigma _R \, ,  \\&
r_{xT} =r_{Tx} =-\frac 2{T} ( \varSigma _T-\frac  \varSigma T)  .
\endaligned \tag 7.3$$
Next notice
$$\frac {2\varSigma _\theta }{T^3 \sinh ^2(R)} \varPsi  _{x\theta}  =
\frac {2\varSigma _\theta }{T^3 \sinh (R)} \sum _{j=1}^2 a_j(Y) \,
\varPsi _{xY_j},$$ with $a_1(Y)   =    -\frac{Y_2}{ \sinh  (R)} $
and $a_2(Y)   =\frac{Y_1}{ \sinh  (R)}  .$ We set for $j=1,2$
$$ r_{xY_j} =r_{Y_jx}= \frac {2\varSigma _\theta }{T^3 \sinh (R)}
a_j(Y). \tag 7.4$$ We then write
$$\aligned &
(P+B^2) \varPsi + \sum _{\alpha } r_{x\alpha }\partial _x
\partial _\alpha   \varPsi   =
-   \frac {P_c}TF_2(\varPsi , \Cal A ) -\\& -\frac
{P_c}{T^2}F_3(\varPsi , \Cal A )  +
    \frac 2T \Cal Q_0(\varSigma ,\Cal A)   {P_c}\phi ^{\prime}    -
  \frac 1T \Cal Q_0(\varSigma ,\varSigma )  {P_c}\th ^{\prime\prime} -\\&
-  \frac 1{T^2}\Cal A\Cal Q_0( \varSigma ,\varSigma ) {P_c} \phi
^{\prime\prime}+  \frac 1T (  \varPsi_{x}+ \Cal A
 \, {P_c}\phi ^{\prime} ) \frac{\Cal H }{\Omega   }
     .
\endaligned \tag 7.5
$$
We now start a long list of inequalities on the terms in the right
hand side of equations (7.1-2) and (7.5). We advise the reader to
skip the remaining part of this section at a first reading and to
come back to these lemmas when they are referenced later.

\proclaim{Lemma 7.2} Assume (6.1) and let $\Cal B_j (\varPsi, \Cal
A)$, $j=2$ (resp. $j=3$) be one of $  G_2(\varPsi, \Cal A)$ and $
H_2(\varPsi, \Cal A)$ (resp. $  G_3(\varPsi, \Cal A)$ and $
H_3(\varPsi, \Cal A)$). Then for $|I|\le N$ we have
$$ \| Z^I \Cal B_j (\varPsi, \Cal A) \| _{L^2_Y}\le
{C( \mu^\prime)\epsilon}   \| (\varPsi , \Cal A)_N \| _{L^2_Y} . $$
\endproclaim
\noindent It is an immediate consequence of Lemma 7.1 and of
Leibnitz rule.

\proclaim{Lemma 7.3} Assume (6.1) and let  $ |I|\le N^\prime $.
Then, for $\Cal B_2 (\varPsi, \Cal A)$  a quadratic expression of
the form $\Cal B_2 = \Cal A \langle \varPsi , \varphi _1 \rangle +
 \langle \varPsi ^2 , \varphi _2 \rangle $ for $
 \varphi _j (x) \in \Cal S (\Bbb R) $,   we
have
$$ \align &
\| Z^I \Cal B_2 (\varPsi, \Cal A) \| _{L^2 _Y} \le
 \frac{C( \mu^\prime)\epsilon ^2} {T^{\frac 12-
2\delta }}   .
\endalign $$
For $\Cal B_3 (\langle \varPsi , \varphi _1 \rangle  , \Cal A) $ a
cubic expression in the arguments, for $ |I|\le N^\prime $ we have
$$ \align &
\| Z^I \Cal B_3 (\langle \varPsi , \varphi _1 \rangle , \Cal A) \|
_{L^2 _Y} \le
  {C( \mu^\prime)\epsilon ^3}     .
\endalign $$
\endproclaim
REMARK. Notice that   $  H_2(\varPsi, \Cal A)$   is by the
discussion after (3.6) of  the above form. Similarly, $G_2(\varPsi,
\Cal A)-3\th \, \Cal A^2 $ is of  the above form.

\noindent The proof of Lemma 7.3 follows from the  Leibnitz rule,
  Lemma 7.1 and (6.1). \qed
\bigskip

\proclaim{Lemma 7.4} Assume (6.1). Then for $|I|+m\le N$ we have
$$ \| Z^I B^m F_j (\varPsi, \Cal A) \| _2\le
{C( \mu^\prime)\epsilon}   \| (\varPsi , \Cal A)_N \| _2 .\tag 1$$
If $ |I|+m\le N^\prime $ we have
$$ \align &
\| Z^IB^m  F_3 (\varPsi, \Cal A) \| _{L_Y^2 W^{1,\frac p{p-1}}_x}
\le C( \mu^\prime) \epsilon \left (  \epsilon +
T^{-\frac{1}{2}+2\delta }\| (\varPsi , \Cal A)_N \| _2\right )
 \| (\varPsi , \Cal A)_N \| _2   .\tag 2
\endalign $$
\endproclaim
 We start with $ Z^I B^m F_2 (\varPsi, \Cal A)$ which is a sum of
terms of the form

$$\align & \big [ B^m C_{J,K}(T,x, Y)\phi^2(x)
\big ] \, (Z^J\Cal{A}(T,Y)) \, (Z^K\Cal{A}(T,Y))
\tag 3
\\ B^m \big  [
&C_{J,K}(T,x,Y)\phi (x)
\, ( Z^J\Cal{A}(T,Y) ) \, (Z^K\varPsi(T,x,Y) ) \big  ]
 \tag 4
\\
&B^m \big  [ C_{J,K}(T,x,Y)\, (Z^J\varPsi(T,x,Y))
\, (
Z^K\varPsi(T,x,Y)) \big  ] .\tag 5
\endalign$$
 with $C_{J,K} \in \Cal E ^0$
and $|J|+|K|\le |I|.$
Now $B^m  \circ C_{J,K} =  C_{J,K} \circ B^m  +[ B^m,C_{J,K}]
  $
with $[ B^m,C_{J,K}]$ a pdo of order $m-1$ by Lemma 2.4. By Lemma 2.4 we have
 $\|
B^m(j)\| _2\le \| (j)\| _{L^2_YH^m_x}$ for $j=3,4,5$  and so  it is
enough to show that for linear combinations $L_j( \varPsi , \Cal A
)= \lambda _j \varPsi (T,x,Y)+\mu _j \psi (x) \Cal A(T ,Y)$, with
bounded coefficients, we have
$$    \| (Z^J\partial _x^j L_1))
\, ( Z^K\partial _x^kL_2) \| _{2}  \le \text{rhs} (1). \tag 6$$   If
say $|K|+k\le \frac N2,$ then by Sobolev embedding we bound by
$$  \| Z^K\partial _x^k L_2\| _\infty
 \| (\varPsi   , \Cal A )_N( T)  \| _{2} \le C
\| L_2\| _{H_Y^{|K|+2} W ^{k+1,p}_x} \| (\varPsi   , \Cal A )_N( T)
 \| _{2}. $$ Since $B$ is an
elliptic pdo, by  $[\frac N2 ] +3\le N^\prime $  and by (6.1)
 $$ \| L_2(T)\| _{H_Y^{|K|+2} W ^{k+1,p}_x} \le C
\| (\lambda _2 \varPsi (T )+\mu _2 \psi (x) \Cal A(T  ) ) _{N^\prime
}\| _{L_Y^2 L^p_x}  \le C \mu ^\prime    \epsilon . \tag 7
$$
hence we have obtained (6).  When we prove (1) for $j=3$ we proceed
similarly reducing to
$$    \| (Z^J\partial _x^j L_1))
\, ( Z^K\partial _x^kL_2) \, ( Z^W\partial _x^wL_3) \| _{2}  \le
\text{rhs} (1) \tag 8$$ with now two factors  with fewer than $\frac
N2 $ derivatives, say $|K|+k\le \frac N2 $ and $|W|+w\le \frac N2.$
For either of them we have a bound like (7) and so we have (8).

 We turn  to the proof of (2).
By Lemma 2.4 and using the above notation,  $$\| B^m  \left (  Z^J
L_1  \,   Z^K L_2  \,  Z^W L_3  \right ) \| _{L^2 _Y W_x^{1,\frac
p{p-1}}} \le C \|  Z^J L_1  \,   Z^K L_2  \,  Z^W L_3\| _{L^2 _Y
W_x^{m+1,\frac p{p-1}}}.$$ Terms of the form $ \psi (x) \Cal
A(\varPsi , \psi (x) \Cal A)^2$ are bounded, by   Schwartz and
Sobolev inequalities and Sobolev embedding and by (6.1), by
$$ \|   (\varPsi , \psi (x) \Cal A) _{ {N^\prime} } \| _{L^2 _Y L_x^ p}
 ^2\| (\varPsi , \Cal A) _N(T)  \| _{2} \le C(\mu ^\prime ) \epsilon ^2\| (\varPsi , \Cal A) _N(T)  \| _{2}. $$
 Terms like $\varPsi ^3$ are bounded  by
 $$  \|    \varPsi   _{ {N^\prime} } \| _{L^2 _Y L_x^ p}  \|  \varPsi  _N(T)  \| _{2}^2
 \le  \mu ^\prime  T^{-\frac{1}{2}+2\delta }\epsilon  \| (\varPsi , \Cal A) _N(T)  \| _{2}^2. \,
 \qed$$
\bigskip
\proclaim{Lemma 7.5}  Let $L=L( \langle \varPsi , \psi \rangle ,
\Cal A )= \lambda \langle \varPsi (T,x,Y) , \psi (x) \rangle
_{L^2_x}+\mu \Cal A(T ,Y)$ with $\lambda$ and $\mu $ two constants,
bounded by a fixed number. Assume (6.1). Then, for $|I|\le N$ we
have $\| Z^I \Cal Q_0( \varSigma , L)\| _2
 \le (7.6)$ with
$$\frac {C( \mu^\prime)\, \epsilon ^2}T+
{C( \mu^\prime)\, \epsilon} \,  (E^{\frac 12} _N(T)+ \tilde E^{\frac
12} _N(T))+\frac {C( \mu^\prime) \epsilon }T  \int _{T_0}^T
(E^{\frac 12} _N(\tau )+ \tilde E^{\frac 12} _N(\tau )) d\tau . \tag
7.6
$$
\endproclaim
\noindent We consider $Z_3^k Z^{I^\prime } \Cal Q_0( \varSigma , L)
,$ where $Z^{I^\prime }$ does not contain $Z_3$. By Lemma 5.3
$$Z^{I^\prime } \Cal Q_0( \varSigma , L)= \sum C_{J^\prime
K^\prime } \Cal Q_0(Z^{J^\prime } \varSigma ,Z^{K^\prime } L)$$ with
$C_{J^\prime K^\prime }$ constants. Suppose $k=0$. If $|J^\prime
|\le [\frac N2]< N^\prime -1,$  by Lemma 7.1 $\| Z^{J^\prime }
\varSigma (T)\| _{L^\infty _Y} \le C(\mu ^\prime ) T^{\frac
12+\delta }
 \epsilon .$
Consequently by (6.1)
$$ \| \Cal Q_0(Z^{J^\prime } \varSigma ,Z^{K^\prime }L)
 (T)\| _{L^2_Y} \le \frac {C(\mu ^\prime ) \epsilon  }{ T^{\frac 12-\delta } }
E^{\frac 12} _N(T).$$ From this point on, we assume   $|J^\prime |
> [\frac N2]$ and  $|K^\prime | \le [\frac N2]$. We claim:
$$ \frac 1{T^2} \| (\partial _R Z^{J^\prime } \varSigma ) \, (
\partial _R Z^{K^\prime } L) +
\frac 1{\sinh ^2(R)}\, (
\partial _\theta  Z^{J^\prime } \varSigma ) \, (
\partial _\theta  Z^{K^\prime } L )\| _2 \le
\frac {C(\mu ^\prime )\epsilon  }{T}
 \tilde E^{\frac 12} (T, \varSigma _N).$$
To show the claim notice that $\| \partial _RZ^{K^\prime } L \|
_\infty \le  \|  (\varPsi , \Cal A) _{[\frac N2]+1} \| _\infty .$
Similarly, equality $\partial _\theta =\tanh (R) \big [- \frac
{Y_2}R Z_1+\frac {Y_1}R Z_2 \big ]$ implies $ \| \frac 1{\sinh
(R)}\,
\partial _\theta Z^{K^\prime } L \| _\infty \le  \|  (\varPsi , \Cal
A) _{[\frac N2]+1} \| _\infty .$ By $N^\prime -2 \ge [\frac N2]+1$
and by Lemma 7.1,  $ \|  (\varPsi , \Cal A) _{[\frac N2]+1} \|
_\infty $ is bounded by (6.1). So our claim holds.

\noindent Still assuming $|J^\prime |
> [\frac N2]$ and $|K^\prime | \le [\frac N2]$, we  consider
$$ \aligned &
 \| (\partial _T Z^{J^\prime } \varSigma - \frac 1T
Z^{J^\prime } \varSigma )
 (\partial _T Z^{K^\prime }L - \frac 1T
Z^{K^\prime }L ) \| _2 \le
 \frac 1{T}\| Z^{J^\prime } \varSigma   \| _2 \|
\partial _TZ^{K^\prime }L \| _\infty \\&  +
 \frac 1{T^2}\| Z^{J^\prime } \varSigma  \| _2   \|
Z^{K^\prime }L  \| _\infty +
 {C(\mu ^\prime ) } \epsilon
\tilde E^{\frac 12} (T, \varSigma _N) ,\endaligned \tag 1
$$
where we expanded in the lhs and used Lemma 7.1. For $T\in
[T_0,T^\ast [$ we have
$$ \aligned  \frac 1T \| \varSigma _N (T) \| _2 \le  & \frac 1T \|
 \varSigma _N (T_0) \| _2+\frac 1T  \int _{T_0}^T
\| \partial _T \varSigma _N (\tau )\| _2 d\tau  . \endaligned \tag
2$$ Since  $|K^\prime | \le [\frac N2],$ by  Lemma 7.1
$$\text{rhs}(1) \le  \frac { C(\mu ^\prime ) \epsilon }T
 (\epsilon +\int _{T_0}^T
\| \partial _T \varSigma _N (\tau )\| _2 d\tau)  +{C(\mu ^\prime ) }
 \epsilon  \tilde E^{\frac 12} (T, \varSigma _N) .$$
So far  we had $k=0$. Let now $k>0$ and set $$ Z_3^k Z^{I^\prime }
\Cal Q_0( \varSigma , L) = \big [ \frac 1T (T\partial _T) \big ] ^k
Z^{I^\prime }\Cal Q_0 (\varSigma , L).$$ By elementary computation
and by Lemma 5.3 this is a sum of terms of the form
$$\frac 1{T^{a} }
\Cal Q_0 ( (T
\partial _T)^{c}Z^{J^\prime } \varSigma ,(T
\partial _T)^{d}Z^{K^\prime } L) $$
where $|J^\prime |+|K^\prime | \le |I^\prime |$, $c+d \le b\le k \le
a$. Distinguishing between  cases $c+ |J^\prime |\le  [\frac N2]$
and $c+ |J^\prime |>  [\frac N2]$, and using the fact that the $T'$s
in the numerator are canceled by the $T'$s in the denominator, we
prove the desired estimate proceeding as in the $k=0$ case.\qed

\bigskip

\proclaim{Lemma 7.6} Using the notation of Lemma 7.5 we have for
$|I|\le N^\prime $
$$\| Z^I \Cal Q_0( \varSigma , L)\| _2
\le
 \frac
{C( \mu^\prime)\epsilon  } {T^{\frac 12-\delta  }}\left ( \epsilon +
\frac{E_N^{\frac{1}{2}}(T)}{T} \right )
 .
 $$
\endproclaim
We consider $Z_3^k Z^{I^\prime } \Cal Q_0( \varSigma , L) ,$ where
$Z^{I^\prime }$ does not contain $Z_3$. By Lemma 5.3
$$Z^{I^\prime } \Cal Q_0( \varSigma , L)= \sum C_{J^\prime
K^\prime } \Cal Q_0(Z^{J^\prime } \varSigma ,Z^{K^\prime } L)$$ with
$C_{J^\prime K^\prime }$ constants. Suppose $k=0$. We claim:
\proclaim{Claim} Consider
$$ \frac 1{T^2} \| (\partial _R Z^{J^\prime } \varSigma ) \, (
\partial _R Z^{K^\prime } L) +
\frac 1{\sinh ^2(R)}\, (
\partial _\theta  Z^{J^\prime } \varSigma ) \, (
\partial _\theta  Z^{K^\prime } L )\| _2 .
  $$ Then $(1) \le \dfrac {C( \mu^\prime)\epsilon  E_N^{\frac{1}{2}}(T)} {T^{\frac
32-\delta }}$ if $|I^\prime |\le N^\prime $.
\endproclaim

 {\it Proof of the Claim.}   $|I^\prime | \le N^\prime
$ implies $|K^\prime | \le N^\prime $ and by Sobolev inequality
$$\|
\partial _RZ^{K^\prime } L \| _\infty  + \| \frac 1{\sinh (R)}\,
\partial _\theta Z^{K^\prime } L \| _\infty \le C\, E_N^{\frac{1}{2}}(T)  . $$
This and the following inequality, consequence of  $|J^\prime | \le
N^\prime $ and (6.1), give us the Claim:
$$\|
\partial _RZ^{J^\prime }  \varSigma \| _2 + \| \frac 1{\sinh (R)}\,
\partial _\theta Z^{J^\prime }  \varSigma \| _2 \le C(\mu ^\prime ) \epsilon  T^{\frac{1}{2}+\delta
} . $$

 We  consider now
$$ \aligned &
 \| (\partial _T Z^{J^\prime } \varSigma - \frac 1T
Z^{J^\prime } \varSigma )
 (\partial _T Z^{K^\prime }L - \frac 1T
Z^{K^\prime }L ) \| _2.\endaligned \tag 1
$$
If $ |   K^\prime |< N^\prime /2$ then by Sobolev embedding and by
(6.1),
$$ \aligned & (1)\le
 \frac 1{T}\| Z^{J^\prime } \varSigma   \| _2 \|
\partial _TZ^{K^\prime }L \| _\infty   + \\& +
 \frac 1{T^2}\| Z^{J^\prime } \varSigma  \| _2   \|
Z^{K^\prime }L  \| _\infty +\|  \partial _T Z^{J^\prime } \varSigma
  \| _2 \|(\partial _T Z^{K^\prime }L - \frac 1T
Z^{K^\prime }L ) \| _\infty \le \\& \le C(\mu ^\prime )\| L
_{N^\prime }\| _2 \left ( \frac{\| Z^{J^\prime } \varSigma   \|
_2}{T}     +   \|
\partial _T Z^{J^\prime } \varSigma
  \| _2    \right ) \le   C(\mu ^\prime )  \epsilon  ^2   T^{-\frac{1}{2}+\delta
} ,\endaligned
$$ where we used (2) in Lemma 7.5, with $N$ replaced by $N^\prime .$

If instead $ |   K^\prime |> N^\prime /2$, then  by (6.1) and by
Lemma 7.1 we have
$$ \aligned & (1)\le
 \frac 1{T}\| Z^{J^\prime } \varSigma   \| _\infty \|
\partial _TZ^{K^\prime }L \| _2   + \\& +
 \frac 1{T^2}\| Z^{J^\prime } \varSigma  \| _\infty   \|
Z^{K^\prime }L  \| _2 +\|  \partial _T Z^{J^\prime } \varSigma
  \| _\infty \|(\partial _T Z^{K^\prime }L - \frac 1T
Z^{K^\prime }L ) \| _2 \le \\& \le C(\mu ^\prime )\epsilon \left (
\frac{\| Z^{J^\prime } \varSigma   \| _2}{T}     +   \|
\partial _T Z^{J^\prime } \varSigma
  \| _2    \right ) \le   C(\mu ^\prime )  \epsilon   ^2 T^{-\frac{1}{2}+\delta
} .\endaligned
$$

Let now $k>0$ and set $$ Z_3^k Z^{I^\prime } \Cal Q_0( \varSigma ,
L) = \big [ \frac 1T (T\partial _T) \big ] ^k Z^{I^\prime }\Cal Q_0
(\varSigma , L).$$ By elementary computation and by Lemma 5.3 this
is a sum of terms of the form
$$\frac 1{T^{a} }
\Cal Q_0 ( (T
\partial _T)^{c}Z^{J^\prime } \varSigma ,(T
\partial _T)^{d}Z^{K^\prime } L) $$
where $|J^\prime |+|K^\prime | \le |I^\prime |$, $c+d \le b\le k \le
a$. Using the fact that the $T'$s in the numerator are canceled by
the $T'$s in the denominator, we prove the desired estimate
proceeding as in the $k=0$ case.\qed

\bigskip
\proclaim{Lemma 7.7} Assume (6.1). Then for $|I|\le N$
 we have
$$\| Z^I  \Cal Q_0( \varSigma ,\varSigma  ) \| _2
\le \frac 1 {T^{\frac 12-\delta }}  (7.6).
$$
\endproclaim

\noindent We write as in Lemma 7.5, $Z^I  \Cal Q_0( \varSigma
,\varSigma  ) =\big [  \frac 1T (T
\partial _T) \big ] ^k
 Z^{I^\prime }  \Cal Q_0( \varSigma ,\varSigma  ) $, where $ Z^{I^\prime } $ does not contain $\partial _T.$
The last expression is a sum of terms of the form
$\frac 1{T^{a}} (T\partial _T) ^b  Z^{I^\prime }  \Cal Q_0,$
 $b+|I^\prime |\le |I|.$ This in turn is a sum of
terms like
$\frac 1{T^{a}} \Cal Q_0( (T\partial _T) ^c Z^J\varSigma ,
 (T\partial _T) ^d Z^K\varSigma  ) $ with $|J|+|K|\le |I^\prime |$
and $c+d\le b$. Since one of $c+|J|$ and $b+|K|$ is less than
$[\frac N2]$ we conclude
$$ \aligned &
\| Z^I  \Cal Q_0( \varSigma ,\varSigma  ) \| _2  \le \left ( \tilde
E ^{\frac 12} (T, \varSigma _{N}) + T^{-1} \|  \varSigma _{N}\| _2
\right ) \times \\& \times
 \|  \, |\partial _T
\varSigma _{ [\frac N2] }|+
T^{-1}| \varSigma _{[\frac N2] }|+T^{-1}|\partial _R \varSigma _{[\frac N2] }|
+ \frac {|\partial _\theta \varSigma _{[\frac N2]}|}{T\sinh ( R)}
) (T) \|_\infty .
\endaligned \tag 1 $$
The second factor in rhs can be bounded by $ \frac {C(
\mu^\prime)\epsilon} {T^{\frac 12-\delta }}$ by
 Lemma 7.1. $T^{-1} \|  \varSigma _N \|_2$ can be bounded
by (2) in Lemma 7.5.\qed
\bigskip
\proclaim{Lemma 7.8} With the notation of Lemma  7.7, we have for
$|I|\le N^\prime $
$$\| Z^I  \Cal Q_0( \varSigma ,\varSigma  ) \| _2
\le
 \frac
{C( \mu^\prime)\epsilon ^2} {T^{1-2\delta }}.
$$
\endproclaim
Indeed we have an inequality like (1) in Lemma
 7.7 with $N^\prime$ replacing $N$.
Then

$$ \|  \left ( |\partial _T
\varSigma _{ [\frac {N^\prime }2] }|+ T^{-1}| \varSigma _{[\frac
{N^\prime }2] }|+T^{-1}|\partial _R \varSigma _{[\frac  {N^\prime
}2] }| + \frac {|\partial _\theta \varSigma _{[\frac  {N^\prime
}2]}|}{T\sinh ( R)} \right ) (T) \|_\infty \le C(\mu ^\prime )
\epsilon T^{-\frac 12+\delta }$$ by Lemma 7.1. We have
$$ \left (  \tilde E ^{\frac 12} (T, \varSigma _{  {N^\prime } }) +
T^{-1}\| \varSigma _ {N^\prime } \| _2 \right ) \le  C(\mu ^\prime )
\epsilon T^{-\frac 12+\delta } $$ by (6.1) and by formula (2) in
Lemma 7.5.\qed

\bigskip

By Leibnitz rule, Lemmas 7.7 and 7.8 and (6.1) imply:
\proclaim{Lemma 7.9} Assume (6.1). Let as in Lemma 7.5, $L=  \lambda
\langle \varPsi (T,x,Y) , \psi (x) \rangle _{L^2_x}+\mu \Cal A(T
,Y)$. Then, for $|I|\le N$
 we have
$$\| Z^I  L \Cal Q_0( \varSigma ,\varSigma  ) \| _2
\le
 \frac
{C(\mu ^\prime )\epsilon } {T^{\frac 12-\delta }}  (7.6).
$$ With the above notation, for  $|I|\le N^\prime $  we have $$
  \| Z^I  L \Cal Q_0( \varSigma ,\varSigma  ) \| _2 \le
 \frac
{C( \mu^\prime)\epsilon ^3} {T^{1-2\delta  }}  .$$
\endproclaim
It is enough to bound for $I=J+K$
$$\| (Z^J  L) Z^K \Cal Q_0( \varSigma ,\varSigma  ) \| _2. \tag 1$$
For $|I|\le N$ for $|J|\le N/2$, by Lemmas 7.1 and 7.7 we have
$$(1) \le \| Z^J  L\| _\infty \|  Z^K \Cal Q_0( \varSigma ,\varSigma  ) \|
_2 \le \frac {C(\mu ^\prime )\epsilon } {T^{\frac 12-\delta }}
(7.6).$$ For $|I|\le N$ for $|J|> N/2$, by Sobolev embedding and by
Lemma 7.8 we have
$$(1) \le \| Z^J  L\| _2\sum _{|H|\le N^\prime } \|  Z^H \Cal Q_0( \varSigma ,\varSigma  ) \|
_2 \le \frac {C(\mu ^\prime )\epsilon ^2} {T^{1-2\delta }}
E^{\frac{1}{2}}_N(T)\le
 \frac
{C(\mu ^\prime )\epsilon } {T^{1-2\delta }}  (7.6).$$ When $|I|\le
N^{\prime}$ for $|J|\le N^{\prime }/2$ by Lemmas 7.1 and 7.8 we have
$$(1) \le \| Z^J  L\| _\infty \|  Z^K \Cal Q_0( \varSigma ,\varSigma  ) \|
_2 \le \frac {C(\mu ^\prime )\epsilon ^3} {T^{1-2\delta }}  .$$ When
$|I|\le N^{\prime}$ for $|J|> N^{\prime }/2$ we have
$$(1) \le \| Z^J  L\| _2\sum _{|H|\le N^\prime } \|  Z^H \Cal Q_0( \varSigma ,\varSigma  ) \|
_2 \le
 \frac
{C(\mu ^\prime )\epsilon  ^3} {T^{1-2\delta }}  .$$

\bigskip
\proclaim{Lemma 7.10} Assume (6.1). Let $|I|+m\le N$. Set
$L=B^m\varPsi, \Cal A , \varSigma .$  Then there is  a fixed
constant $C$such that
$$\| [Z^I,P] L \| _2 \le \frac C{T^2} \left (    E^{\frac 12} _N(T)+      \tilde
E^{\frac 12} _N(T)    \right ).
$$
\endproclaim
We write $[Z^I,P] L = [\partial _T^k,P] Z^{I^\prime } L $ where $
Z^{I^\prime }$ does not contain $\partial _T$. This is a sum of
terms of the form
$$  \frac 1{T^{k_1+2} } \Delta _{hyp } Z^{\tilde I }
L , \quad  Z^{\tilde I }=\partial _T^{k_2} Z^{I^\prime } \, , \tag
1$$ with $k_1 >0$,  $k_1 +k_2=k$ and $ |\tilde I |  \le |I|-1.$ Then
apply to the terms in (1) the following lemma:

\proclaim{Lemma 7.11} For $L$ as in Lemma 7.10 we have for $|I|+m\le
q-1$ and for a fixed $C$
$$\| \frac 1{T  } \Delta _{hyp } Z^{  I }
L \| _2\le C \left ( E^{\frac 12} _q(T)+      \tilde E^{\frac 12}
_q(T) \right )
$$
where the $L^2$ norm is either the $L^2_{xY}$ one for $L= B^m\varPsi
$ or the $L^2_{ Y}$ one for $L=   \Cal A , \varSigma .$
\endproclaim

Write $\| \frac 1{T  } \Delta _{hyp } Z^{  I } L  \| _2=\| \dots \|
_{R\le 1}+\| \dots  \| _{R\ge 1}$. Then $$\| \frac 1{T  } \Delta
_{hyp } Z^{  I } L  \| _{R\ge 1}\lesssim E^{\frac 12} _q(T)+ \tilde
E^{\frac 12} _q(T) $$ by $\Delta _{hyp}=
\partial ^2_R +\frac 1{ \tanh (R)}
\partial _R +\frac 1{ \sinh  ^2(R)}  \partial ^2_\theta $.
 To prove
$\| \frac 1{T  } \Delta _{hyp } Z^{  I } L  \| _{R\le 1}\lesssim
E^{\frac 12} _q(T)+      \tilde E^{\frac 12} _q(T) $, by (5.2) we
write
$$ \frac 1{T }  |\Delta _{hyp }
 Z^{  I } L |\le
| \frac {Z_1}T Z_1
 Z^{  I } L |
+
| \frac {Z_2}T Z_2
 Z^{ I } L |
+
| \frac {Z_0}T Z_0
 Z^{  I } L |$$
and see that by the formulas
$$\aligned & Z_1=\frac {Y_1}R \partial _R -
\frac 1{\tanh (R)} \frac {Y_2}R\partial _\theta  \quad , \quad
Z_0 =\partial _\theta \\&
 Z_2=\frac {Y_2}R \partial _R +
\frac 1{\tanh (R)} \frac {Y_1}R\partial _\theta
\endaligned
$$
  we have that $
\frac 1{T }  \|\Delta _{hyp }
 Z^{  I }L \| _{L^2(R\le 1) }$ is bounded by a fixed constant times
$$
\aligned & \| \sum _{j=1}^2| \frac {\partial _R}T  Z_j
 Z^{I  }L |+\sum _{j=0}^2
| \frac {\partial _\theta } {T \sinh (R)} Z_j Z^{I  } L | \| _2 \le
C \left ( E^{\frac 12} _q(T)+      \tilde E^{\frac 12} _q(T) \right
) .\qed
\endaligned
$$

\bigskip
\head \S 8 Energy inequalities \endhead We set
$$\aligned E_{1\varSigma } (T, \varPsi _q)=&  E_{{(1)} } (T, \varPsi _q) +
\sum _{m+|I| \le q} \int \big \{ 2r_{Tx } \partial _x (B^m
Z^{I}\varPsi )
\partial _T (B^m Z^{I}\varPsi )-\\&  -\sum _{\alpha , \beta }r_{\alpha \beta }
 \partial _\alpha (B^m Z^{I}\varPsi )
  \partial _\beta  (B^m Z^{I}\varPsi ) \big \}
\sinh (R) \, dxdRd\theta ,
\endaligned $$
with $\alpha $ and $\beta $ summed over all $T,R, x$ and $Y_j$,
$j=1,2$, with $r_{\alpha \beta }=0$ if it is not in the   list
   (7.3-4). We set $E_{\varSigma }(T,(\varPsi,\Cal{A})_q)=
E_{1\varSigma }+ E_{(2)}$. For $\epsilon _0(\mu ^\prime )$ small,
(6.1) and Lemma 7.1
  imply for any $q\le N$
$$ \frac 12 E_q(T) \le E_{\varSigma }(T,\varPsi _q) \le 2E_q(T) \quad \forall
\, T\in [T_0,T^\ast [ \quad \forall \, q\le N. \tag 8.1$$

 \proclaim
{Lemma 8.1} There is an $\epsilon (\mu ^\prime ) >0 $ and   fixed
constants $C_0$ and $C_1$ such that if (6.1) is true for $\epsilon
\in ]0, \epsilon (\mu ^\prime ) [$
 then for any $T\in [T_0,T^\ast [$ we have:
$$\aligned & E(T,(\varPsi, \Cal A )_N)+\tilde E(T, \varSigma _N) \le
C_0 T^{C_1C(\mu ^\prime ) \epsilon } \epsilon ^2.
\endaligned $$
\endproclaim
\noindent REMARK. By adjusting $C(\mu ^\prime )$ we simply write
$E_N+\tilde E_N \le C_0 T^{C(\mu ^\prime ) \epsilon } \epsilon ^2.$
\bigskip
\noindent PROOF of Lemma 8.1. Let us set $$ \aligned & \frac {C(
\mu^\prime)\epsilon  ^2} {T} +\left ( \frac {C( \mu^\prime)\epsilon
} {T}+\frac {C} {T^2}\right )
 \big [ E _N^{\frac 12}
+\tilde E _N^{\frac 12} \big ] (T) +\frac {C( \mu^\prime) \epsilon }
{T^  {\frac{5}{2}-\delta }} E _N (T)  \\&  +  \frac {C( \mu^\prime)
\epsilon } {T^2}\int _{T_0}^T \big [ E _N^{\frac 12} +\tilde E
_N^{\frac 12} \big ](\tau ) d\tau .
\endaligned \tag 8.2$$ Then we have:

\proclaim{Lemma 8.2} Given (6.1)   in  $[T_0,T^\ast [$ for any $q$,
$q\le N$, we have $\frac d{dT}\tilde E^{\frac 12}(T,\varSigma _q)
\le (8.2) $.
\endproclaim

\proclaim{Lemma 8.3}  Given (6.1)   in  $[T_0,T^\ast [$ for any
  $q$, $q\le N$ we have $
 \frac d{dT}E^{\frac 12}_{\varSigma }(T,(\varPsi,\Cal{A})_q) \le
 (8.2) .
$
\endproclaim
We assume Lemmas 8.2 and 8.3 and continue the proof of Lemma 8.1. By
a continuity argument we assume initially that the last term in the
first line of (8.2), that is $\frac {C( \mu^\prime) \epsilon } {T^
{\frac{5}{2}-\delta }} E _N (T)$, is not present. Setting $D= E
^{\frac 12}_\varSigma (T, ( \varPsi , \Cal A )_N)+ \tilde E ^{\frac
12}(T, \varSigma _N )$, $\varphi (T)=2(\frac {C(
\mu^\prime)\epsilon} {T}+\frac {C}{T^2})$, adding from Lemmas 8.2
and 8.3 and using (8.1) we obtain
$$\frac d{dT}D \le  \frac
{2C( \mu^\prime)\epsilon  ^2} {T} +\varphi (T) D+  \frac {2 C(
\mu^\prime)\epsilon } {T^2}\int _{T_0}^T D(\tau ) d\tau .$$
Integrating,
$$D(T)\le D(T_{0})+ 2C( \mu^\prime)\epsilon  ^2\log \frac{T}{T_0}+\int _{T_0}^T \varphi (\tau )D(\tau )d\tau +
\int _{T_0}^T \frac {2 C( \mu^\prime)\epsilon } {\tau }D(\tau )d\tau
.
$$
By Gronwall inequality we obtain an inequality of the desired form
$$D(T)\le e^{ \int _{T_0}^T \left ( \frac {4 C( \mu^\prime)\epsilon
} {\tau } + \frac{C}{\tau ^2}\right ) d\tau }\left ( D(T_{0})+ 2C(
\mu^\prime)\epsilon ^2\log \frac{T}{T_0} \right ) .$$ By a
continuity argument it is easy now to absorb the $\frac {C(
\mu^\prime) \epsilon } {T^ {\frac{5}{2}-\delta }} E _N (T)$ term
inside $\frac{C}{T^2}E _N^{\frac{1}{2}} (T).$ \qed

\bigskip

We return now to Lemmas 8.2 and 8.3. Apply $Z^I$ to (7.1)
$$\align PZ^I\varSigma
=&   Z^I  \frac{\Cal H}{\Omega } +[Z^I,P] \varSigma .\tag 8.3
\endalign
$$
Taking $m\ge 0$ with $B^0=P_c$, apply $B^mZ^I$ to (7.5):
$$\aligned &
(P+B^2)B^mZ^I\varPsi + \sum _{\alpha } r_{x\alpha }\partial _x
\partial _\alpha B^mZ^I \varPsi   =
-B^mZ^I\big \{ \frac 1TF_2(\varPsi , \Cal A ) +\\& +\frac
1{T^2}F_3(\varPsi , \Cal A )\big \} +
 Z^I\big [ \frac 2T \Cal Q_0(\varSigma ,\Cal A)\big ] B^m\phi ^{\prime}    -
Z^I \big [ \frac 1T \Cal Q_0(\varSigma ,\varSigma ) \big ]
 B^m \th ^{\prime\prime} -\\&
- Z^I\big [\frac 1{T^2}\Cal A\Cal Q_0( \varSigma ,\varSigma )\big
]B^m \phi  ^{\prime\prime}+Z^IB^m\big [ \frac 1T (  \varPsi_{x}+
\Cal A
 \, \phi ^{\prime}) \frac {\Cal H}\Omega  \big ] \\&
-  [Z^I,P]B^m \varPsi  - \sum _{\alpha} [Z^I,r_{x\alpha }] B^m
\varPsi _{x\alpha } - \sum _{\alpha } r_{x\alpha }[B^m,\partial _x
\partial _\alpha  ] Z^I \varPsi .
\endaligned \tag 8.4
$$ Finally there is a similar equation for $\Cal A$.

$$\aligned  &(P+\frac 32) Z^I\Cal A =\\&
= - Z^I\left [ \frac 1TG_2(\varPsi , \Cal A)+\frac 1{T^2}G_3(\varPsi
, \Cal A)+ \frac 1{T}\Cal Q_0( \varSigma ,\varSigma )  + \frac 1T
\langle \varPsi ,\phi  ^{\prime \prime } \rangle \frac {\Cal
H}\Omega \right ]  \\& -Z^I \left [ \frac 2T \Cal Q_0( \varSigma, L(
\varPsi ,\Cal A ) )  + \frac 1{T^2}\Cal Q_0( \varSigma ,\varSigma
)L( \varPsi ,\Cal A )\right ] -[Z^I,P] \Cal A.
\endaligned \tag 8.5
$$

By Lemma 6.1, Lemmas 8.3 and  8.4 follow from:

\proclaim{Lemma 8.4} For $j=3,4,5$ we have $\| rhs (8.j)\| _2\le
(8.2).$
\endproclaim
Lemmas 7.2 and 7.4 guarantee, for $m+|I|\le N$,  that for $\ell
=2,3$, quadratic or cubic expressions succinctly denoted by  $ (
\varPsi , \Cal A )^\ell $ satisfy
$$\| B^mZ^I
 \frac 1 {T^{\ell -1}  }  ( \varPsi , \Cal A )^\ell \| _ 2 \le \frac{C( \mu^\prime)}{T^{\ell -1}}\epsilon E
_N^{\frac 12}(T).\tag 1$$ Lemma  7.3   guarantees for $ |I|\le
N^\prime $
$$\|  Z^I
 \frac 1 {T^{\ell -1}  }  ( \langle \varPsi ,
 \psi \rangle , \Cal A )^\ell \| _ 2 \le \frac{C( \mu^\prime)}{T^{\ell -1}}\epsilon  ^\ell .\tag 2$$
Lemma 7.5 guarantees   for $ |I|\le N  $
$$\|  Z^I\frac{1}{ T }\Cal Q_0(
\varSigma , L(  \langle \varPsi ,
 \psi \rangle  ,\Cal A ))\| _ 2 \le \frac{(7.6)}{T}\le
(8.2)\tag 3$$ and Lemma 7.6 guarantees for $ |I|\le N ^\prime  $
$$\aligned & \|  Z^I\frac{1}{ T }\Cal Q_0(
\varSigma , L(  \langle \varPsi ,
 \psi \rangle ,\Cal A  )\| _ 2 \le \frac{C(\mu ^\prime ) \epsilon  }{T^{\frac{3}{2} -\delta }}
 \left ( \epsilon +\frac{E_N^{\frac{1}{2}}(T)}{T} \right ) \text{ for } |I|= N ^\prime\\&
\|  Z^I\frac{1}{ T }\Cal Q_0( \varSigma , L(  \langle \varPsi ,
 \psi \rangle ,\Cal A  )\| _ 2 \le \frac{C(\mu ^\prime ) \epsilon ^2 }{T^{\frac{3}{2} -\delta }}
  \text{ for } |I|< N ^\prime
\endaligned
 .\tag 4$$
Lemma 7.7 guarantees for $ |I|\le N  $
$$\|  Z^I\frac{1}{ T }\Cal Q_0(
\varSigma , \varSigma )\| _ 2 \le \frac{(7.6)}{T^{\frac{3}{2}
-\delta }}\le (8.2) \tag 5$$ and Lemma 7.8 guarantees for $ |I|\le
N^\prime  $
$$\|  Z^I\frac{1}{ T }\Cal Q_0(
\varSigma , \varSigma )\| _ 2 \le  \frac{C(\mu ^\prime ) \epsilon
^2}{T^{2 -2\delta }}. \tag 6$$ Lemma 7.9 guarantees for $ |I|\le N
   $
$$\|  Z^I\frac{1}{ T^2 }L(  \langle \varPsi ,
 \psi \rangle ,\Cal A )\Cal Q_0(
\varSigma , \varSigma )\| _ 2 \le \frac{(7.6)}{T^{\frac{5}{2}
-\delta }}\le (8.2)\tag 7$$ and for $ |I|\le N^\prime
   $
$$\|  Z^I\frac{1}{ T^2 }L(  \langle \varPsi ,
 \psi \rangle ,\Cal A )\Cal Q_0(
\varSigma , \varSigma )\| _ 2 \le  \frac{C(\mu ^\prime ) \epsilon
^3}{T^{3 -2\delta }}.  \tag 8$$ The following Lemma holds:

\proclaim{Lemma 8.5} We have
$$\align &\| (\Cal H )_N\| _ 2\le  (8.2)\, , \quad  \| (\Cal H ) _{N^{\prime }}\| _ 2\le  \frac{C(\mu ^\prime ) \epsilon
 }{T^{\frac{3}{2} -2\delta }} \left ( \frac{E_N^{\frac{1}{2}}(T)}{T}+\epsilon \right )  ;\tag 8.6 \\& \| \left ( {\Cal
H}/{\Omega } \right )_N\| _ 2\le  (8.2)\, , \quad  \| \left ( \Cal H
/{\Omega } \right )  _{N^{\prime }}\| _ 2\le  \frac{C(\mu ^\prime )
\epsilon
 }{T^{\frac{3}{2} -2\delta }} \left ( \frac{E_N^{\frac{1}{2}}(T)}{T}+\epsilon \right ) .\tag 8.7\endalign $$ We
have
$$ \| \left ( \frac {\Cal H  -\frac{6}{T}\Cal A \langle \varPsi , \th \th
^\prime \phi \rangle }  {\Omega } \right )  _{N^{\prime }}\| _ 2\le
\frac{C(\mu ^\prime ) \epsilon
 }{T^{\frac{3}{2} - \delta }} \left ( \frac{E_N^{\frac{1}{2}}(T)}{T^{1-\delta }}+\epsilon \right )  .\tag
8.8$$ We have
$$ \| \left ( \text{rhs}(8.5)+Z^I\frac{3}{T}\Cal A ^2 \| \th \phi ^3
\|_1  +[Z^I,P]\Cal A\right )  _{N^{\prime }}\| _ 2\le \frac{C(\mu
^\prime ) \epsilon  }{T^{\frac{3}{2} -2\delta }} \left (
\frac{E_N^{\frac{1}{2}}(T)}{T}+\epsilon \right ) .\tag 8.9$$
\endproclaim
(8.6) follows from inequalities from (1) to (8) and from Lemma 7.3.
(8.7) follows from (8.6) and Leibnitz rule. (8.9) follows from
inequalities from (1) to (8),  Lemma 7.3 and   (8.7).  (8.8)
follows from (1) to (8), from Lemma 7.3 and Leibnitz rule because we
have eliminated the only term which decays
 like $T^{-\frac{3}{2}
+2\delta }$, the $  \Cal A \langle \varPsi ,
  \th \th
^\prime \phi \rangle  $ one.  (8.10) is proved like the previous
ones, exploiting the  fact that $N_1<N^{\prime}$, which allows to
exploit the inequalities with $<N^\prime $ derivatives in Lemma 7.6
and in inequality (4) in Lemma 8.4.

 \bigskip

We resume the proof of Lemma 8.4.
 Lemma 7.10 guarantees for $m+ |I|\le N
   $ $$\| [Z^I,P]   (
B^m\varPsi ,\Cal A , \varSigma ) \| _ 2 \le \frac{C}{T^{2}}\big ( E
_N^{\frac 12} +\tilde E _N^{\frac 12} \big )\le (8.2). \tag 9$$
Hence (9) and (8.7) guarantee Lemma 8.4 for $j=3$.  (8.7) and the
Leibnitz rule guarantee for $|I|\le N$
$$ \| Z^I \left ( \frac{1}{T} \langle \varPsi , \psi \rangle \frac{\Cal H}{\Omega
} \right ) \|  \le  (8.2). \tag  10$$ The above estimates guarantee
Lemma 8.4 for $j=5$.

\noindent Given a linear combination $L=\lambda \varPsi (x,Y,T) +
\mu \Cal A ( Y,T) \psi (x)$, we want to show that for $|I|+m\le N$
we have
$$\| Z^I B^m \partial _x L \frac{\Cal H}{T\Omega }   \|_2\le (8.2) .
\tag 11$$ We consider for $|J|+|K|\le |I|$
$$\| (Z^J B^m \partial _x L) (Z^K \frac{\Cal H}{ T\Omega })  \|_2. \tag 12$$
If $|K| \le [ N ]/2 $ then  the desired inequality follows from
$$\aligned & (12) \le
  \| Z^J B^m
\partial _x L \| _2 \| Z^K \frac{\Cal H}{T\Omega } \|_\infty \le  C
T^{-1} E_N ^{\frac{1}{2}} (T) \|  (  {\Cal H}/{ \Omega })_{N^\prime
} \|_2\\& \le \frac{C(\mu ^\prime ) \epsilon  }{T^{\frac{5}{2}
-2\delta }}E_N ^{\frac{1}{2}} (T)    \left (
\frac{E_N^{\frac{1}{2}}(T)}{T}+\epsilon \right )       \le (8.2)
  .\endaligned $$
If $|K| \ge  [ N ]/2+1 $ then $|J|+m \le  [N/2]-2.$ Then $$  (12)
\le \| Z^J B^m
\partial _x L \| _\infty \| Z^K \frac{\Cal H}{T\Omega } \|_2 \le  C
\| Z^K \frac{\Cal H}{T\Omega } \|_2 \sup _{|J|+m \le N^\prime -4}
\| Z^J
\partial_x ^m L \| _2 .$$
Then $\| Z^K \frac{\Cal H}{T\Omega } \|_2 \le    T^{-1}(8.2)$ by
(8.7) and $ \| Z^J
\partial_x ^m L \| _2\le C T^{\frac 12 - \frac 1p } \| Z^J
\partial_x ^m L \|  _{L^2_y L^p _x}$  by a lemma   stated
and proved immediately below, see Lemma  8.6. These last two
estimates give $(12)\le (8.2)$ also for $|K| \ge  N^\prime $. We
state and prove Lemma 8.6 and then we continue the proof of Lemma
8.4.

\proclaim {Lemma 8.6 } Let $1\le  p\le q \le \infty $   and fix any
$j$ and $ J .$   Assume (6.1). Then there are constants $C$ and
$C_M$ such that for any $T\in [T_0,T^\ast [$, any $Y$  we have for
our solution: $\|
\partial _x^jZ^J \varPsi (T,Y) \| _{L^p _x} \le $
$$
\le C T^{\frac 1p - \frac 1q } \left ( \| \partial _x^jZ^J \varPsi
(T,Y) \| _{L^q _x} + \frac {C_M}{T^M} (|Z^J\varSigma (T,Y)|+|Z^J\Cal
A (T,Y)|) \right )
 .$$
 \endproclaim
For simplicity let us pick $j+|J|=0$, but the general   argument is
the same. First of all we write $\int_{\Bbb R} dx|
\varPsi(T,x,Y)|^{p} =\int_{|x+\sigma (t,y)|\le T}+\int_{|x+\sigma
(t,y)|\ge T}.$ Then, by H\"older,
$$\int_{|x+\sigma (t,y)|\le T}dx| \varPsi(T,x,Y)|^{p}
\le T^{1 - \frac pq } \| \varPsi(T,x,Y) \|_{L^q_x}^{p}.$$  Next we
consider the $|x+\sigma (t,y)|\ge T$ integral. In intervals of
existence, the solution $w$ of (1.4) is supported in the set
$x^2+y^2=x^2+T^2\sinh ^2 R \le (t+K)^2 < (t+2K)^2=T^2\cosh ^2 R .$
This implies that $T \ge |x|$ on the support of $w$. By (1.10) and
(3.2)
$$\psi(t,x,y)=w(t,x+\sigma (t,y),y)+\sigma (t,y) \int_0^1
\th ^\prime (x+\sigma (t,y)-s \sigma (t,y))ds-\phi (x) a(t,y).$$ So,
for  $|x+\sigma (t,y)|\ge T,$
$$\psi(t,x,y)=\sigma (t,y) \int_0^1\th ^\prime
(x+s \sigma (t,y))ds-\phi (x) a(t,y).$$ Lemma 7.1 and the equality
$\sigma =\frac \varSigma T$ imply (here we are focusing on the
$\sigma$'s in the argument of $w$ and  $\th ^\prime$)
$$| \sigma  (t,y)|\le
C \mu ^\prime  T^{-\frac 12+\delta } \epsilon .$$ Therefore, for
$\epsilon $ small,  inequality $|x+\sigma (t,y)|\ge T$ implies
$|x|>T/2$. Therefore, for some fixed constants $C $ and $C_M$,
$$\aligned &\int_{|x+\sigma (t,y)|\ge T}dx| \varPsi(T,x,Y)|^p \le  \\& \le
 C|\varSigma (T,Y)| ^p  \int_{|x|\ge T/2}dx  \, |\th ^\prime (x)|^p
+ C|\Cal A (T,Y)|^p \int_{|x|\ge T/2}dx| \phi (x)|^p
\\&\le C_MT^{-M} \big (|\varSigma (T,Y)|^p
+ |\Cal A (T,Y)| ^p \big  ) ,
\endaligned$$
with $M>0$ an arbitrarily large number.\qed
\bigskip
To complete the proof of Lemma 8.4 we still need to prove $\|  rhs
(8.4)\| _2 \le (8.2).$ We have $\| B^mZ^I  \frac 1TF_2(\varPsi ,
\Cal A )\|_ 2\le (8.2) $ by (1), $\| B^mZ^I  \frac 1{T^2}F_3(\varPsi
, \Cal A )\|_ 2\le T^{-1}(8.2) $ by (1), $\| B^mZ^I  \frac 1{T }\Cal
Q_0(\varSigma ,\Cal A) \phi ^{\prime}  \|_ 2\le  (8.2) $ by (3), $\|
B^mZ^I  \frac 1{T }\Cal Q_0(\varSigma ,\varSigma ) \th
^{\prime\prime} \|_ 2\le  (8.2) $ by (5), $\| B^mZ^I \frac
1{T^2}\Cal A\Cal Q_0( \varSigma ,\varSigma )  \phi ^{\prime\prime}
\|_ 2\le T^{-\frac{3}{2}+\delta } (8.2) $ by (7), $\| [Z^I,P]B^m
\varPsi  \|_ 2\le (8.2) $ by (9) and $\| B^mZ^I \big [ \frac 1T (
\varPsi_{x}+ \Cal A
 \, \phi ^{\prime}) \frac {\Cal H}\Omega
\|_ 2\le   (8.2) $ by (12). Hence  we have $\|  rhs (8.4)\| _2 \le
(8.2) $ if we prove:

\proclaim {Lemma 8.7} Assume (6.1). Then we have for $|I|+m\le N$:

$$ \align &  \| [Z^I,r_{x\alpha } ] B^m \varPsi _{x\alpha }
\| _2 \le  (8.2)\tag 1\\&
  \| r_{x\alpha }[\partial _x \partial _\alpha ,B^m ] Z^I \varPsi \| _2
\le  (8.2)\tag 2
 \endalign $$
\endproclaim
Recall $$\aligned & r_{xx} =\frac 1{T^2} \Cal Q_0( \varSigma ,
\varSigma ) \quad ,
  \quad
r_{xR}=r_{Rx} =\frac 2{T^3}  \varSigma _R \, ,  \\& r_{xT} =r_{Tx}
=-\frac 2{T} ( \varSigma _T-\frac  \varSigma T)  \quad , r_{xY_j}
=r_{Y_jx}= \frac {2\varSigma _\theta }{T^3 \sinh (R)} a_j(Y)
  \quad
\endaligned $$
   where
$a_1(Y)   =    -\frac{Y_2}{ \sinh  (R)} $ and $a_2(Y)   =\frac{Y_1}{
\sinh  (R)}  .$  We prove (1) first. $ [Z^I,r_{x\alpha } ] B^m
\varPsi _{x\alpha }$
 is a sum of terms  $(Z^JB^m \varPsi _{x\alpha })(
 Z^Kr_{x\alpha }),$ with $m+|J|+|K|\le |I|.$
If $m+|J| \le  [\frac N2]$, by Lemma 7.1, $N^\prime \ge [\frac N2]
+3$,  Lemma 8.6 and (6.1) we have

$$\| Z^JB^m \varPsi _{x\alpha }\| _
{L^\infty_YL^p_x} \le C(\mu ^\prime )  ( T^{\frac 12-\frac 1p} \|
Z^JB^m \varPsi _{x\alpha }\| _ {L^\infty_YL^2_x} +\frac \epsilon T)
\le C(\mu ^\prime ) \epsilon T^{-\frac 1p+2\delta }.$$  We have $ \|
Z^K r_{xx } \| _2 \le  T^{-1} (8.2) $ by (5) in  Lemma 8.4, $C(\mu
^\prime ) \epsilon \| Z^K r_{x\alpha  } \| _2  \le  T^{-1} (8.2) $
for $\alpha = R,Y_j$ by the definitions, and $ C(\mu ^\prime )
\epsilon  \| Z^K r_{xT  } \| _2 \le  T^{-1} (8.2) $ which follows by
$ \| Z^K \varSigma _T/T \| _2 \le T^{-1}\tilde E_N^{\frac{1}{2}} (T)
$ and by formula (2) in Lemma 7.5.

\noindent If $m+|J| >  [\frac N2]$, then $ \| Z^K r_{xx } \| _2 \le
C(\mu ^\prime ) T^{-3+2\delta } $ by (6). $  \| Z^K r_{x\alpha } \|
_2 \le   \mu ^\prime   T^{-2 } \epsilon $ for $\alpha \neq x$ by
(6.1). We have $\| Z^JB^m \varPsi _{x\alpha  }\| _2\le CE^{\frac 12}
(T,\varPsi _N), $ where we use $|J|+m\le N-1$ and $\| Z^JB^m \varPsi
_{x\alpha }\| _2\le C\| Z^JB^{m +1 }\varPsi _{\alpha  }\| _2$ by
$\varPsi=P_c(H)\varPsi .$  All these estimates imply claim (1) of
Lemma 8.7.

\noindent  We prove (2) in Lemma 8.7.    For $\alpha =T,R, Y_j$,
  since
$[\partial _x, B^m]$ is a pdo of order $m$,  since $B$ is elliptic
and $\varPsi =P_c(H) \varPsi$
$$ \|  r_{x\alpha } [\partial _x, B^m]
\partial _\alpha Z^I \varPsi \| _2  \le C  \|  r_{x\alpha }\| _\infty
\| B^m \partial _\alpha Z^I \varPsi \| _2 \le \frac {C(\mu ^\prime
)\epsilon } {T^{\frac{3}{2} -\delta }} E^{\frac 12}_N.$$
 For $\alpha = x$,  since $[\partial
_x^2, B^m]$ is a pdo of order $m+1$,
$$ \|  r_{xx } [\partial _x^2, B^m] Z^I \varPsi \| _2  \le C
 \|  r_{xx }\| _\infty \| B^{m+1} Z^I \varPsi \| _2 \le
\frac {C(\mu ^\prime )\epsilon ^2} {T ^{3-2\delta }} E^{\frac
12}_N.\qed
$$

\bigskip
 \head \S 9 The elimination of the $Y$ variable \endhead

Following Klainerman, see \cite{K,DFX}, the energy estimates in the
previous section are used now to interpret terms in the equations
for $\varPsi$ and $\Cal A$, that is (7.2) and (7.3), with
derivatives in $Y$, as perturbations integrable in $T$. Hence the
equation for $\varPsi$ is interpreted as a Schr\"odinger equation
with time $T$ and space variable $x$, the equation for $\Cal A$ is
interpreted as an ODE with time $T$. Specifically, we write
$$
\varPsi_{TT}
+B^2\varPsi +r_{xT} P_c\varPsi _{Tx} = -\frac{1}{T}P_c F_2(\varPsi,\Cal{A})
-\frac{1}{T^2}R_\varPsi
,
 \tag 9.1 $$
with
$$\aligned R_\varPsi &=P_c F_3(\varPsi,\Cal{A})
+\Delta_{hyp}\varPsi +T\Cal Q_0(\varSigma ,   \varSigma  ) P_c \th
^{\prime \prime } -\\& -
 2T \Cal Q_0(\varSigma ,\Cal A) P_c\phi  ^\prime  +\Cal A\Cal
Q_0( \varSigma ,\varSigma )P_c  \phi   ^{\prime \prime }+\\& +
Q_0( \varSigma ,\varSigma )  P_c\varPsi _{xx}
+
T P_c(  \varPsi_{x}+ \Cal A
 \, \phi  ^\prime )P\varSigma  -T^2
\sum _{\alpha \neq T}   r_{x\alpha } \varPsi_{x\alpha }
\endaligned
$$
and we write
$$
\Cal{A}_{TT}
+\frac 32\Cal{A}=-\frac{3}{T}\Cal{A}^2
\| \th \, \phi^3\|_1
-\frac{1}{T^2}R_\Cal{A},
 \tag 9.2 $$
with

$$\aligned
R_
\Cal{A}&=  T\langle 3\th \,\varPsi^2
+6\th \,
\Cal{A}\, \phi\,\varPsi,\phi \rangle +
G_3(\varPsi,\Cal{A})
 + \Cal Q_0(\varSigma ,   \varSigma  )
\langle \varPsi , \phi ^{\prime \prime }\rangle +T\langle \varPsi ,
\phi ^{\prime \prime }
\rangle P\varSigma +\\&
+ 2T\Cal Q_0( \varSigma ,\langle \varPsi , \phi  ^\prime \rangle  )
+ \Cal A Q_0( \varSigma ,\varSigma  )\|\phi  ^\prime \| ^2_2
+ \Delta_{hyp}\Cal{A} +T\Cal Q_0(\varSigma ,\varSigma )
\| \th
^{\prime\prime } \phi \|_1
.
\endaligned
$$

\noindent We have:

\proclaim{Lemma 9.1} There are a fixed constant $C$ and an $\epsilon
_0>0$ such that for $|I|\le N^\prime $ and for $0<\epsilon <\epsilon
_0$ we have for all $T\in [T_0, T^\ast [$:
$$\aligned
\| Z^I R_\Cal{A}(T) \|_{L^2_Y} \le  C T^{\frac 34 } \epsilon
.\endaligned  $$

\endproclaim
Lemma 9.1 is consequence of (8.9) and, by (5.2) and Lemma 8.1, for
$L=\Cal A , \varPsi $ and for $|I| \le N-2$, of
$$\| Z^I\Delta _{hyp}L \| _2 \le \|  L_{|I|+2} \| _2\le \|
L_{N} \| _2\le C T^{C(\mu ^\prime ) \epsilon }\epsilon .\tag 9.3$$

\bigskip The following two lemmas are proved in \S 14.
  \proclaim{Lemma 9.2} There are a fixed constant $C$ and an $\epsilon
_0>0$ such that   for $0<\epsilon <\epsilon _0$ and for $|I|+m\le
N^\prime $ we have  for all $T\in [T_0, T^\ast [$
$$\| Z^IB^{m+1} R _\varPsi  \| _{L^2 _{xY}} \le C T^{\frac 34}
\epsilon .$$
\endproclaim \bigskip \proclaim{Lemma 9.3} There are a fixed constant $C$ and an $\epsilon
_0>0$ such that   for $0<\epsilon <\epsilon _0$ and for $|I|+m\le
N^\prime $ we have   for all $T\in [T_0, T^\ast [$:
$$
 \| B^m Z^I R_\varPsi (T)
\|_{L^2_Y W^{1,\frac p{p-1}}_x} \le C T^{\frac 34} \epsilon .\qed
$$
\endproclaim

\head \S 10 Normal form argument for $\Cal A $ \endhead
\bigskip
The starting points are (9.2) and Lemma 9.2. Consider $\Cal{A}_{\pm}
$ related to $\Cal{A}$ by

$$\Cal{A}_{\pm}= (\mp i \partial_T+\sqrt{3/2})\Cal{A}
\quad, \quad
\Cal{A}=\frac{\Cal{A}_{+}+\Cal{A}_{-}}{\sqrt{6}}.$$
Then  write
$$\aligned &(\pm i \partial_T+\sqrt{3/2})\Cal{A}_{\pm}
= -3\frac{\| \th \, \phi^3\|_1} {6\, T } (\Cal{A}_{+}+\Cal{A}_{-})^2
 -\frac{1}{T^2} R_\Cal{A} .
\endaligned $$
By the theory of normal forms there are constants $\alpha _\pm $,
$\beta _\pm $, $\gamma _\pm $
so that, if we write
$\Cal{A}_{\pm }=A_{1+} +\frac {\alpha _\pm }T( A_{1\pm})^2+
\frac { \beta _\pm }T  A_{1\pm}A_{1\mp }+
\frac { \gamma _\pm }T  (A_{1\mp })^2$,
 we obtain

$$\aligned &(\pm i \partial_T+\sqrt{3/2}) A_{1,\pm}
= -\frac{1}{T^2} R_{\Cal{A}\pm},
\endaligned $$
where we have, by  Lemma 9.1:

\proclaim{Lemma 10.1} There is a fixed constant $C$ such that for
$R_\Cal{A}$ as above we have for all $T\in [T_0, T^\ast [$:
$$\sum_{|I|\le N^\prime } \| Z^I R_{\Cal{A}\pm}(T)  \|_{L^2_Y}
\le  C T^{\frac 34} \epsilon .\qed $$
\endproclaim
\noindent   We have $  R_{\Cal{A}+}-  R_{\Cal{A}}= O(
 A_{1,\pm}^3)+O( T^{-1}A_{1,\pm}R_{\Cal{A}} ) $ plus smaller
terms.  By Leibnitz rule,  Sobolev   inequality,( 6.1), Lemma 8.1
and $\epsilon$ small
  $\| Z^I (A_{1,+},A_{1,-})^3\| _2\le $ $$
  \| (A_{1,+},A_{1,-})_N\| _2^2\| (A_{1,+},A_{1,-}) _{N^\prime }\|
_2\approx \| (\Cal A_{ +} ,\Cal A_{ -})_N\| _2^2\| (\Cal A_{ +}
,\Cal A_{ -}) _{N^\prime }\| _2\le  C  T^{C(\mu ^\prime )\epsilon}
\epsilon  .$$ We have by (6.1), Lemma 9.1 and $\epsilon$ small
$$\| Z^I T^{-1}
(  R_{\Cal{A}}A_{1,\pm } ) \| _2\lesssim  T^{-1}\| (\Cal A_{ \pm } )
_{N^\prime }\| _2  \|  (  R_{\Cal{A}}) _{N^\prime }\| _2\le C
\epsilon  T^{-\frac{1}{4}}.$$

\bigskip
We have: \proclaim{Lemma 10.2} Assume (6.1). Then there is a fixed
constant $C$ such that $\forall$ $I$ with $|I|\le N^\prime $ we have
for all $T\in [T_0, T^\ast [$
$$ \| Z^I A_{1\pm }(T)\|_{L ^{2}_Y}\le C
 \epsilon .$$
\endproclaim

For any multiindex $I$ with $|I|\le N^\prime $ we
have
$$(\pm i\partial_T+\sqrt{3/2})Z^I A_{1\pm}= \frac{1}{T} F^I_{\pm}
\qquad{\text{with}}\quad
\| F^I_{\pm} (T) \|_{L^2 _Y} +
\| F^I_{\pm} (T) \|_{L^\infty  _Y}\le C\epsilon T^{-\frac 14}.$$
Next, we can write
$$\frac 12 \partial_T | Z^I A_{1\pm}|^2 \le \frac{1}{T}
\Im \big \{F^I_{\pm} Z^I\bar A_{1\pm}\big \}$$
from which we obtain
$$ \frac d{dT} \|  Z^I A_{1\pm} \| _{L^2_Y} \le
 \frac{1}{T}\| F^I_{\pm} (T) \|_{L^2 _Y} \le  C\epsilon T^{-\frac 54}. \qed$$

\bigskip

 \head
\S 11 Bounds  on $\Psi$ \endhead

  We have:

\proclaim{Lemma 11.1} Assume (6.1) and Lemmas 9.2 and 9.3.
 Then there are a fixed $c$ and an
 $\epsilon _0(\mu ^\prime ) $  such that if
$\epsilon \in ]0,  \epsilon _0(\mu ^\prime )[ $ and $T \in [T_0,
T^\ast [$ and for any choice of multiindex  $I$ and nonnegative
integer $m$ with $|I|+m\le N^\prime $
$$
\| B^m Z^I\varPsi (T) \|_{L^2 _YL^p_x} \le c  T^{ 2\delta -\frac
12} \epsilon .
$$
\endproclaim
We start from (9.1). Apply $ B^mZ^I$ and consider
$$(B^mZ^I\varPsi)_{TT}+B^2B^mZ^I\varPsi = -Z^I \frac {B^m}T F_2(\varPsi ,\Cal A)-Z^I
 \frac {B^m}{T^2}  R_\varPsi -Z^I r_{xT}B^m \varPsi
_{xT} . $$ We  write  $B^mZ^I\varPsi (T) =
\sum_{j=1}^4B^m\varPsi_{Ij} (T)$ where
$$\aligned
&B^m\varPsi_{I1} (T)=\frac{\sin(B(T-T_0))} B B^mZ^I\varPsi _T(T_0)+
\cos (B(T-T_0)) Z^IB^m\varPsi (T_0)
\\
&B^m\varPsi_{Ij} (T) =-\int_{T_0}^T \frac{\sin B(T-\tau)} {B}\Cal
F_j (\tau)\, {d\tau}\, , \quad j=2, 3 ,4 \\& \Cal F_2(T) = {P_c}Z^I
\frac 1T B^m F_2 (\Cal{A} ,\varPsi ), \quad \Cal F_3(T)=  Z^I \frac
1{T^2} B^m R_\varPsi (T) \\&
 \Cal F_4(T)= Z^Ir_{xT}  B^m\varPsi _{Tx} (T)  .
\endaligned$$
Then we have: \proclaim {Lemma 11.2} We have:
$$\align &
\| B^m \varPsi_{I1} (T) \|_{L^2 _YL^p_x} \le c \langle T
\rangle^{-\frac 12+\frac 1p} \epsilon  \tag 1 \\& \| B^m
\varPsi_{I2} (T) \|_{L^2 _YL^p_x} \le C(\mu ^\prime ) \langle T
\rangle^{-\frac 12+\frac 1p} \log (\langle T \rangle) \epsilon ^2
 \tag 2 \\&
\| B^m \varPsi_{I3} (T) \|_ {L^2 _YL^p_x} \le C(\mu ^\prime )
\langle T \rangle ^{-\frac 12+\frac 1p} \epsilon
 \tag 3 \\&
\| B^m \varPsi_{I4} (T) \|_ {L^2 _YL^p_x} \le C(\mu ^\prime )
\langle T \rangle ^{-\frac 12 +2\delta } \epsilon ^2.
 \tag 4
\endalign
$$\endproclaim
\noindent  (1) follows from
  Proposition 4.1 and Corollary 2.3.  (3) Follows from
$$\align &\int_{T_0}^T
\| \frac{\sin B(T-\tau)} {B}B^mZ^I R_\varPsi (\tau )\|_{L^2
_YL^p_x} \frac{d\tau}{\tau ^2} \le \\& \le  \int_{T_0}^T \| B^mZ^I
R_\varPsi (\tau )\| _{L_Y^2 W^{1,\frac p{p-1}}_x}
\frac{Cd\tau}{(T-\tau )^{\frac 12-\frac 1p}\tau ^2} \le
\int_{T_0}^T\frac{ C_1 \epsilon d\tau}{(T-\tau )^{\frac 12-\frac
1p}} \tau ^{-\frac {10}9},
\endalign$$
by Lemma 9.3. As for  inequality (2), notice that if we consider
$$ \int_{T_0}^T
\| \frac{\sin B(T-\tau)} {B} P_cZ^I \frac 1\tau B^mF_2 (\Cal{A}
\phi,\varPsi )\|_ {L^2 _YL^p_x} {d\tau}, \tag 5
$$
then by Corollary 2.3,
$$
(5) \le c\int_{T_0}^T ( T-\tau)^{-\frac 12+\frac 1p} \| P_cZ^I
\frac 1\tau B^mF_2 (\Cal{A} \phi,\varPsi )\|_{L^2_YW_x^{1,\frac
p{p-1}}} {d\tau} .$$ Now $\tau Z^I \frac 1\tau B^mF_2$
 is formed by terms
schematically of the form
$$\align &  \left [ Z^J\Cal{A}(\tau,Y)\right  ] \left [ Z^K\Cal{A}(\tau,Y) \right  ]  B^m
\left [ C_{J,K}(\tau,x,Y)\phi^2(x) \right  ]
  \tag 6
\\
& \left [ Z^J\Cal A  (\tau,Y)\right ]   B^m \left [
C_{J,K}(\tau,x,Y)\phi (x) Z^K \varPsi(\tau,x,Y)\right ]   \tag 7
\\
& B^m \left [C_{J,K}(\tau,x,Y)Z^J\varPsi(\tau,x,Y)
Z^K\varPsi(\tau,x,Y)\right  ] .\tag 8
\endalign $$ with $C_{J,K} \in \Cal E ^0.$ In (6) $|J|+|K|\le N^\prime$.
By (6.1), Lemma 7.1 and Lemma 2.4 we have
$$\aligned
&\| (6) \|_{L^2_YW_x^{1,\frac p{p-1}}} \le \|  C_{J,K}(\tau,x,Y)
\phi ^2 \| _{L^\infty_YW_x^{m+1,\frac p{p-1}}} (\mu ^\prime )^2
\epsilon ^2.
\endaligned$$ Similarly
  $ \| (8) \|_{L^2_YW_x^{1,\frac p{p-1}}}$
is bounded by a sum , for $j+k \le m+1$,   for $|K|+k\le N/2$ and
for $\frac{1}{r}= \frac  {p-1}p-\frac{1}{p}$, and by Sobolev
embedding in the second step,

$$\aligned & C\,
\| Z^JB^j\varPsi(\tau ) \| _{L^\infty_YL_x^{r}} \|
Z^KB^k\varPsi(\tau ) \| _{L^2_YL_x^{p}} \le  E _{(1)} ^{\frac{1}{2}}
(\tau , \varPsi )  \, \| Z^KB^k\varPsi(\tau ) \| _{L^2_YL_x^{p}} .
\endaligned  $$ By (6.1) and  Lemma  8.1, the latter $\le C(\mu ^\prime ) \epsilon
^2\tau ^{2C(\mu ^\prime ) \epsilon  -(\frac 12-2\delta ) } $, which
is bounded. A similar bound is obtained for (7), and thus we obtain
(2).
\bigskip

Now we turn to the proof of (4) Lemma 11.2.  $\Cal F _4$ is a sum of
terms of the form

$$\frac 1{T^{1+a}}\, ( \partial _T^b Z^J \varSigma ) \,
(Z^K B^m \varPsi _{Tx} ) \tag 9 $$
 with $|J|+|K|\le  |I|,$ $b\le
1$ and $a\ge 1$ if $b=0$. If $|J|\le [\frac N2]$ then by Lemma 7.1
$$\| Z^J \frac 1T (\varSigma _T - \frac \varSigma T ) \| _\infty
\le C(\mu ^\prime ) \epsilon T^{-\frac 32+\delta }$$ while by Lemmas
2.4 and  8.6 and by $|K|+m+2\le N^\prime +2\le N-1$
$$\|  Z^K B^m \varPsi _{Tx}\| _{L^2 _Y W^{1,\frac p{p-1}}_x} \le C
\tau ^{\frac 12 -\frac 1p} ( \|  \varPsi _{N^\prime +2} \| _{L^2 _Y
H^1_x} +\frac {C(\mu ^\prime ) \epsilon }{\tau ^2})\le C \tau
^{\frac 12 -\frac 1p+C(\mu ^\prime ) \epsilon }C(\mu ^\prime )
\epsilon .$$
 As a consequence the desired estimate  follows from
the upper bound
$$ \int_{T_0}^T \frac{ (9)}{(T-\tau)^{\frac 12-\frac 1p}}\le
C(\mu ^\prime ) \epsilon \int_{T_0}^T
 \frac 1{(T-\tau)^{\frac 12-\frac 1p}}
 \frac{d\tau}{\tau ^{1+\frac 1p -2\delta}}.\tag 10$$
If  in (9)  $|J|> [\frac N2]$, then by (6.1)
$$\| Z^J \frac 1T (\varSigma _T - \frac \varSigma T ) \| _2
\le C(\mu ^\prime ) \epsilon T^{-\frac 32+\delta }$$ while by Lemma
8.6 and Sobolev embedding
$$\|  Z^K B^m \varPsi _{Tx}\| _{L^\infty  _Y W^{1,\frac p{p-1}}_x} \le C
\tau ^{\frac 12 -\frac 1p} ( \|  \varPsi _{N} \| _{L^2 _Y
L^2 _x} +\frac {C(\mu ^\prime ) \epsilon }{\tau ^2})\le C \tau
^{\frac 12 -\frac 1p+C(\mu ^\prime ) \epsilon }C(\mu ^\prime )
\epsilon .$$ We obtain again (10). \qed
\bigskip

Along with the $L^2_YL^p_x$ estimate, in (6.1) we have also low
energy estimates for $\varPsi. $ We have:

\proclaim{Lemma 11.3}Assume (6.1). There is a fixed      constant $C
$ and a constant
 $\epsilon _0(\mu ^\prime ) $  such that if
$\epsilon \in ]0,  \epsilon _0(\mu ^\prime )[ $ and $T \in [T_0,
T^\ast [$ we have
$$  \| \partial _T \varPsi _{N^\prime} (T) \|
_{L^2 _Y H^1_x}+ \|  \varPsi _{N^\prime} (T) \| _{L^2 _Y H^2_x}\le 2
\epsilon $$
\endproclaim

  We apply  $B^mZ^I$ to formula (9.1)  for
$m+|I|\le N^\prime +1$ and $m\ge 1$:
$$ \aligned &
(B^{m }Z^I\varPsi)_{TT}+B^2B^{m }Z^I\varPsi + r_{xT}
\partial _xB^{m}Z^I\varPsi _T=\\& =
-Z^I \frac 1T B^{m} F_2(\varPsi ,\Cal A)-Z^I
 \frac 1{T^2} B^{m} R_\varPsi -\\& -[  Z^I ,  r_{xT}  ]   B^{m}
\varPsi _{xT} - r_{xT}[ B^{m}   ,\partial _x   ] Z^I\varPsi _T.
\endaligned  \tag 11.1$$
We have: \proclaim{Claim} We have $$\aligned & \| \text{rhs}
(11.1)+Z^I \frac 1T B^{m }  3\th \, \phi ^2(x) \Cal A ^2  \| _{L^2
 _{xY} } \le \frac{C(\mu ^\prime )}{T^{\frac 98}} \epsilon ^2 .
\endaligned
$$
\endproclaim
Assume the Claim.    For $D(T)= \| ( B^{m}Z^I \varPsi )_T \| ^2
_{L^2 _Y L^2_x }+ \| B^{m+1}Z^I \varPsi \| ^2 _{L^2 _Y L^2_x }$ by
  the Claim   we get for $\psi (x)=3\th \, \phi ^2(x)$

$$\aligned &\frac 12 \frac{d}{dT} D(T)  \le   D^{\frac 12}(T)
\frac{ C(\mu ^\prime )\epsilon ^2} {T^{\frac 98}}  - \langle (
B^{m}Z^I \varPsi )_T, Z^I \frac 1T B^{m }  \psi   \Cal A ^2  \rangle
_{L^2  _{xY}} .
\endaligned $$
If on an interval $[T_0,T_1[$ we have $D(T)\le 8 \epsilon ^2$,
 then for $\epsilon \in ]0, \epsilon (\mu ^\prime )[ $
we have $$ \aligned & D(T)\le D(T_0)+  C(\mu ^\prime )\epsilon
^3\int _{T_0}^T \frac{ d\tau } {\tau ^{\frac 98}} -     \frac 1\tau
 \langle  B^{m}Z^I \varPsi  (\tau ) ,B^m \psi (x)
  Z^I \Cal A ^2(\tau )\rangle  _{L^2  _{xY}} \mid  _{\tau =T_0}^{\tau =T}  \\& + \int
  _{T_0}^T
 \langle  B^{m}Z^I \varPsi  (\tau ) ,B^m \psi (x)
  \partial _TZ^I \frac 1\tau\Cal A ^2(\tau )\rangle  _{L^2  _{xY}} .\endaligned \tag 1$$

Now we claim
$$
 \big | \langle  B^{m}Z^I \varPsi  (T ) ,B^m \psi (x)
  \partial _TZ^I \frac 1T\Cal A ^2(T )\rangle  _{L^2  _{xY}}\big |
  \le
\frac{ C(\mu ^\prime )\epsilon ^3} {T^{\frac 54}} .\tag 2
  $$
By (6.1)   we have $\| \partial _TZ^I \frac 1T\Cal A ^2(T )\| _2\le
\frac{C(\mu ^\prime )\epsilon ^2}{T}.$   (2) will follow from
  $\|  B^{m}Z^I \varPsi   \|  _{L_Y^2 L_x^p }\le
C(\mu ^\prime ) T^{-\frac 14} \epsilon .$ Recall $m+|I|\le N^\prime
+1$ and $m\ge 1$. By interpolation $$\aligned  &\|    B \varPsi
_{N^\prime }\| _{ L_x^p } \le C_1\|     \varPsi _{N^\prime }\|  _{
L_x^p }^{\frac{N-N^\prime -2}{N-N^\prime -1}} \|       \varPsi _{N-
 1}\|  _{  L_x^p }^{\frac{1}{N-N^\prime -1}}
\\& \le C_2\| \varPsi  _{N^\prime }\| _{  L_x^p }^{\frac{N-N^\prime -2}{N-N^\prime
-1}} \|
  \varPsi _{N   }\|  _{  L_x^2
}^{\frac{1}{N-N^\prime -1}}\le  C(\mu ^\prime )\epsilon  \,
T^{\frac{ C(\mu ^\prime )\epsilon}{N-N^\prime -1} - \frac{N-N^\prime
-2}{N-N^\prime -1 } (\frac{1}{2}-2\delta ) }   \\& \le C(\mu ^\prime
)\epsilon    \, T^{ {\tilde C(\mu ^\prime )\epsilon} - \frac{6}{7 }
(\frac{1}{2}-2\delta ) }.
\endaligned \tag 3
$$ where for the second inequality we
use Sobolev embedding and for the second (6.1) and Lemma 8.1. From
(3) we get (2). Entering the information in  (1), we get $D(T)\le
\epsilon ^2+      {C}(\mu ^\prime ) \epsilon ^2 (
T^{-\frac{1}{8}}_{0} + \epsilon )$ for some fixed function $ {C}(\mu
^\prime ).$ Since $T_0$ can be thought large and $\epsilon >0$
small, we conclude that $D(T)\le 4 \epsilon ^2$ for $T\in
[T_0,T_1]$. So we  have proved that for any  $T_1<T^\ast ,$ $D(T)\le
8  \epsilon ^2$   in $  [T_0,T_1]$ implies $D(T)\le 4  \epsilon ^2$
in $  [T_0,T_1]$. Hence we conclude $D(T)\le 4 \epsilon ^2  $ in $
[T_0,T^\ast [$.

\bigskip We prove the claimed inequality (2).
By Sobolev Embedding, by (7.3) and by $m+|I|\le N^\prime +1$ with
$m\ge 1$ for the first inequality, and by (6.1) and Lemma 8.1 for
the second,
$$\aligned &
\| [  Z^I ,  r_{xT}  ]   B^{m} \varPsi _{xT}
 \| _{2}
\le \frac 1T(\| \partial _T \varSigma _{N^\prime  } \| _2+ \frac 1T
\| \varSigma _{N^\prime  } \| _2) \| \varPsi \| _{H^N} \le C(\mu
^\prime ) \epsilon ^2 T^{C(\mu ^\prime )\epsilon   +\delta
-\frac{1}{2} }.
\endaligned $$
The rhs is bounded by $C(\mu ^\prime ) T^{-\frac 98 } \epsilon ^2$.

\noindent By Lemma 7.1, by the fact that $[  B^{m},
\partial _x ]$ is a pseudodifferential operator of order $m$
and by Lemma 8.1 ,
$$\frac 1T
\| (\varSigma _T- \frac \varSigma T) \| _{\infty} \| [ B^{m},
\partial _x ]  Z^I\varPsi _T \| _{2 } \le \frac{
C(\mu ^\prime )\epsilon }{T^{1+\frac 12-\delta }}
 \| \varPsi \| _{H^N}\le  \frac{ C(\mu ^\prime )\epsilon ^2 }
{T^{1+\frac 12-\delta - C(\mu ^\prime )\epsilon }}.
$$
By Lemma 9.2, which we have yet to prove, $  \| Z^{I}\frac 1{T^2}
B^{m} R_\varPsi \| _2 \le C(\mu ^\prime ) \epsilon ^2
T^{-\frac{9}{8}}$.

\noindent By $ F_2(\varPsi ,\Cal A)-3\th \, \phi ^2(x) \Cal A ^2= 6
\th \, \phi (x)\Cal A \varPsi + 3\th (x)\, \varPsi ^2,$  we claim
$$\aligned &
\| Z^I \frac 1T B^{m } \big [ F_2(\varPsi ,\Cal A)-3\th \, \phi
^2(x) \Cal A ^2 \big ] \| _{2 } \le \frac{C(\mu ^\prime )}{T^{\frac
98}} \epsilon ^2 .
\endaligned \tag 4
$$
To check (4) observe that we need to bound a combination of $\varPsi
^2$ and of $\Cal A \varPsi $.  We have $\| ( \varPsi ^2)_{N^\prime
+1}\| _2 \le C \| \varPsi
 _{N^\prime  }\|  _{L^2 _Y L^p_x}\|   \varPsi
 _{N  }\|  _{2}\le   C(\mu ^\prime ) T^{-\frac 12 +2\delta +C(\mu ^\prime )\epsilon } \epsilon ^2 $
by (6.1), H\"older inequality, Sobolev embedding and Lemma    8.1.
Next, we recall that $m+|I|\le N^\prime +1$ with $m\ge 1$. Hence it
is enough to bound
$$\| B( \psi (x)\Cal A \varPsi )_{N^\prime  }\| _2 \le C \|  \Cal A _{N }\|
_2  \|  (1+B) \varPsi _{N^\prime   }\| _{L^2_YL^p_x}       .\tag 5
$$
By Lemma 8.1 we have $ \|  \Cal A _{N }\| _2 \le   C(\mu ^\prime )
T^{ C(\mu ^\prime )\epsilon } \epsilon$. By (3) and by (6.1) we have
$$ \| (1+B) \varPsi _{N^\prime   }\| _{ L^2_YL^p_x}  \le    \|   \varPsi _{N^\prime   }\| _{ L^2_YL^p_x} +
 \| B \varPsi _{N^\prime   }\| _{ L^2_YL^p_x}  \le  C(\mu
^\prime )\epsilon   T^{ {\tilde C(\mu ^\prime )\epsilon} -
\frac{6}{7 } (\frac{1}{2}-2\delta ) } .$$
\bigskip

\bigskip
\head \S 12  Estimates for $\varSigma$
\endhead
\noindent In Klainerman's classical proof of dispersion of solutions
of the zero mass wave equation, see texts \cite{Ho,So}, the Morawetz
vectorfield takes a central role. Notice though that in these
treatements dimension is at least 3. With this remark in mind, we
recall for $|I|\le N^\prime $ equation $PZ^I\varSigma = rhs(8.3).$
We then consider
$$\sinh (R) \, \, \Cal K(Z^I\varSigma ) \, \,
PZ^I\varSigma =\Cal K(Z^I\varSigma  )\, \,
 rhs(8.3) $$
and set, in the notation of Lemma 5.4,
$$\Cal E_1(T,Z^I\varSigma  )= \iint dR \, d\theta \,
\frac{ \Cal P^0}{\cosh (R)}.$$ Then we have:

\proclaim{Lemma 12.1} Assume (6.1). Then, for $\epsilon _0$ small
enough and for $\epsilon \in ]0, \epsilon _0[$, for $|I|\le N^\prime
$ there is a fixed $C$ such that
$$\aligned &\frac d{dT}\Cal E_1(T,Z^I\varSigma) \le ( C(\mu ^\prime )
\epsilon + C ) \epsilon \, T^{-\frac 12+\delta  }   \Cal E ^{\frac
12} _1(T,Z^I\varSigma ) +\frac{1} {T} \Cal E_1(T,Z^I\varSigma ).
\endaligned$$
\endproclaim
Let us assume for the moment Lemma 12.1. By Gronwall inequality and
by (6.1) we get:

\proclaim {Lemma 12.2} Assume (6.1) Then, there are a fixed $C$ and
 an $\epsilon _0$ small enough such that for $\epsilon \in ]0, \epsilon
_0[$, for $|I|\le N^\prime $ we have
$$ \Cal E_1^ {\frac 12}(T,Z^I\varSigma)\le CT^{\frac 12+\delta} \epsilon ^.$$
In particular, for $|I|\le N^\prime $,
$$\align & \| \frac 1{ \sinh  (R)} \partial _\theta Z^I\varSigma \| _{L^2 _Y}
\le C T^{\frac 12+\delta} \epsilon \tag 1 \\& \| T\partial
_TZ^I\varSigma + \tanh (R)\,  \partial _RZ^I\varSigma \| _{L^2 _Y}
\le  C T^{\frac 12+\delta} \epsilon .\tag 2 \endalign $$
\endproclaim

\bigskip
Notice that (1) can be used to get (6.2) for $\frac 1{T \sinh  (R)}
\partial _\theta  \varSigma$ but that (2) is not enough for $\partial
_T \varSigma$ and for $ \frac{\partial _R }{T}\varSigma $. For this
reason we introduce
$$\Cal
E_2(T,Z^I\varSigma ) =\frac 12 \iint d\theta \, dR\, \sinh  (R) \big
[ T^2(\partial _TZ^I\varSigma )^2+(\partial _RZ^I\varSigma )^2 +
\frac{(\partial _\theta Z^I\varSigma ) ^2}{\sinh ^2 (R)} \big ] .$$
Then we claim:

\proclaim {Lemma 12.3} We have for a fixed $C_1$, $\frac d{dT}\Cal
E_2(T,Z^I\varSigma )\le$ $$\le ( C(\mu ^\prime ) \epsilon + C_1
)\epsilon \,  T^{-\frac 12+\delta  }  \Cal E_2^{\frac
12}(T,Z^I\varSigma )  +\frac 1T \iint d\theta \, dR\, \sinh (R) \,
(T\partial _TZ^I\varSigma )^2 .$$
\endproclaim
\bigskip We postpone the proof.
As a consequence of Lemmas 12.2  and 12.3 we obtain: \proclaim
{Lemma 12.4} For $|I|\le N^\prime$ and a fixed constant $\tilde A$,
we have
$$
\align & \| T\partial _T Z^I\varSigma  \| _{L^2_Y}\le \tilde A  T
^{\frac 12+\delta } \epsilon \tag i \\& \| \partial _R Z^I\varSigma
\| _{L^2_Y}\le \tilde A T ^{\frac 12+\delta }\epsilon . \tag ii
\endalign
$$
\endproclaim
\bigskip
PROOF Notice that $\Cal E_2$ is defined with a factor $1/2$, so the
proof does not follow immediately from Lemma 12.3 and Gronwall
inequality, so we use  Lemma 12.2. We   express $[T_0, T^\ast [ $ as
a   union of intervals $[T_1, T_2[ \subseteq [T_0, T^\ast [ $ such
that we have one of the following alternatives: \roster
\item   $\forall \, T\in [T_1, T_2[$ we have
$ \| \partial _R Z^I\varSigma  \| _{L^2_Y}^2 \le \frac{1-2\delta}{
1+ 2\delta  } \| T\partial _T Z^I\varSigma  \| _{L^2_Y}^2$;

\item  $\forall \, T\in [T_1, T_2[$ we have
$ \| \partial _R Z^I\varSigma  \| _{L^2_Y}^2 \ge \frac{1-2\delta}{
1+ 2\delta  } \| T\partial _T Z^I\varSigma  \| _{L^2_Y}^2$.
\endroster

\noindent The union can be taken maximal, in the sense that
intervals of type (1) and (2) alternate.
 For
$[T_1, T_2[ $ as in  (1), from (2) in Lemma 12.3 we conclude the
  that Lemma 12.4 holds  in $[T_1,
T_2[ $ for  $\tilde
 A= \tilde
 B\doteqdot {C}\left (   1-\sqrt{ 1-2\delta}/\sqrt{
1+ 2\delta  }   \right ) ^{-1}$, $C$ of Lemma 12.2 and $\tilde B$
defined by the equality. For   $[T_1, T_2[ $  as in  (2), then
$$   \| T\partial _T Z^I\varSigma  \| _{2}^2 \le       \frac{1 +2 \delta }{2}
(\| T\partial _T Z^I\varSigma  \| _{2}^2+\| \partial _R Z^I\varSigma
\| _2^2).
$$

By Lemma 12.3,
 $\forall \, T\in [T_1, T_2[$ we have:
$$\frac d{dT}\Cal E_2(T,Z^I\varSigma )\le ( C(\mu ^\prime )
\epsilon + C_1 )\epsilon \,  T^{-\frac 12+\delta  } \Cal E_2^{\frac
12}(T,Z^I\varSigma ) +\frac  {1}{ T} \left ( 1 + 2\delta \right
)\Cal E_2(T,Z^I\varSigma )$$ and $ \frac d{dT}\Cal E_2^{\frac
12}(T,Z^I\varSigma )\le \frac 12 ( C(\mu ^\prime ) \epsilon + C_1
)\epsilon \,  T^{-\frac 12+\delta  } +    \frac{1 + 2\delta }{2T}
\Cal E_2^{\frac 12}(T,Z^I\varSigma ).
 $
By Gronwall we get $ \Cal E_2^{\frac 12}(T,Z^I\varSigma ) \le (\frac
T {T_1})^{\frac 12+\delta }\Cal E_2^{\frac 12}(T_1,Z^I\varSigma ) +(
C(\mu ^\prime ) \epsilon ^2+ C_1 \epsilon )
 T^{\frac 1{2}+\delta } $. Since $T_1$ is the endpoint of an
interval of type (1), we have   $\Cal E_2^{\frac
12}(T_1,Z^I\varSigma )\le  \tilde B T^{\frac 12+\delta }_1
\epsilon$. Then $ \Cal E_2^{\frac 12}(T,Z^I\varSigma ) \le
 (\frac T {T_1})^{\frac 1{2}+\delta } \tilde B{T_1}^{\frac 12+\delta } \epsilon
+( C(\mu ^\prime ) \epsilon ^2+ C_1 \epsilon )  T^{\frac 1{2}+\delta
}  $. For $\epsilon $ small,   the latter gives Lemma 12.4 with
$\tilde A = \tilde B +2C_2 $.
\bigskip
What is left now is the proof of Lemmas 12.2 and 12.3. By Lemma 5.5
we have
$$\aligned \frac d{dT}\Cal E_1(T,Z^I\varSigma) &\le
\| (T\partial _T+\tanh R \partial _R) Z^I\varSigma \| _{L^2_Y} T\|
\text{rhs} (8.3) \| _{L^2 _Y}
 -\\&  -\frac{1} {T} \iint dR \, d\theta \frac{ T \tanh (R)} {\cosh (R)}
\Cal P^1 .\endaligned $$ By Lemma 5.6 we have, since the divergence
terms in $\partial _R$ and $\partial _\theta$ disappear and by
definition of $\Cal E_2$,
$$\frac d{dT}\Cal E_2(T,Z^I\varSigma )\le  2\Cal E_2^{\frac 12}(T,Z^I\varSigma
)  T\| \text{rhs} (8.3) \| _{L^2 _Y}+\frac 1T \iint d\theta \, dR\,
\sinh (R) \, (T\partial _TZ^I\varSigma )^2 .$$ We have now:

 \proclaim{Lemma 12.5} Let $|I|\le  N^\prime .$ Then we have
$$ \| \text{rhs}(8.3)\| _2 \le ( C(\mu ^\prime )
\epsilon + C )\epsilon \,   T^{-\frac 32+\delta }.$$
\endproclaim

Notice that Lemma 12.3 follows immediately from Lemma 12.5, while
Lemma 12.2  is the consequence of the following two claims, after
whose proof we start the proof of Lemma 12.5.  \proclaim{ Claim 1}
We have
$$-T \frac{ \tanh (R)} {\cosh (R)}\Cal P^1 \le \frac{\Cal P^0}{\cosh (R)} .$$
\endproclaim
Indeed
$$\aligned  &- \frac{ \Cal P^1} {\cosh (R)}
 \le  \frac 1{2T}\sinh (R)   \times \\& \times
\left [  \tanh (R) (T\partial _TZ^I\varSigma )^2+
 2T\partial _TZ^I\varSigma \,  \partial _RZ^I\varSigma +
 \tanh R (
\partial _RZ^I\varSigma )^2
\right ] \endaligned
$$
and $$\aligned  &\frac {\Cal P^0}{\cosh (R)}
\ge  \frac 1{2} \tanh (R) \cosh (R)\times \\& \times
\left [   (T\partial _TZ^I\varSigma )^2+
 2T\tanh (R)\, \partial _TZ^I\varSigma \,  \partial _RZ^I\varSigma +
  (
\partial _RZ^I\varSigma )^2
\right ] .\endaligned
$$
From the definition of $\Cal P^0$ we obtain:
\proclaim{ Claim 2}  We have $$
\sinh (R)  \, (T\partial _TZ^I\varSigma + \tanh (R)
 \partial _RZ^I\varSigma )^2\le 2\Cal P^0/\cosh (R) .
$$
\endproclaim

\noindent PROOF of Lemma 12.5. First of all, Lemmas 7.10 and 8.1
give us $\| [Z^I ,P]\varSigma \| _2 \le C T^{-2+C(\mu ^\prime )
\epsilon } \epsilon .$ By (8.8) in Lemma 8.5, Lemma 12.5  will hold
if   $  \| Z^I  \frac{\Cal A }{T}\langle \varPsi, \psi \rangle  \|
_2  \le \frac {C(\mu ^\prime ) \epsilon ^2}{T^{\frac 32 - \delta }}
$. The last inequality    is a consequence of the following
refinement of (8.6), which gives us a $\delta $   gain in the decay
of $\varPsi$: for $|I|\le N^\prime $ and $\psi (x)$ rapidly
decreasing, we have
  $$  \| Z^I  \langle \varPsi, \psi  \rangle  \| _2  \le
  {C(\mu ^\prime ) \epsilon  }{T^{-\frac 12 + \delta }}. \tag
12.1$$ (12.1) is crucial in our argument. Indeed  when we estimate
$\varSigma$ the exponents are tight. A decay $T^{-\frac 12 + 2\delta
}$ in (12.1), would lead to a decay $T^{-\frac 12 + 2\delta }$ for
$\varSigma$, with a disastrous feedback  effect on the other
estimates, also on (12.1). So let us prove (12.1) assuming (6.1).
Using the notation in \S 11 we write
$$Z^I  \langle \varPsi, \psi  \rangle =  \sum _{j=1}^4
\langle Z^I\varPsi _{ j},    \psi \rangle .$$ The terms
corresponding to  $j\neq 4$ satisfy (12.1) by   Lemma 11.2. Indeed
it is for   $j=4$  that we need to gain a $T^{-\delta }$.   We have
$$
  \langle Z^I\varPsi _{ 4},  \psi  \rangle  =-2
\int _{T_0}^T Z^I \left [\langle \varPsi _{Tx}(\tau ) , \frac {\sin
(T-\tau )B}B
  \psi \rangle (\varSigma _T -\frac \varSigma \tau )
\frac {1 }\tau \right ] d\tau  .$$ Ignoring the $-2,$ this   is a
sum of terms of the form
$$ \int _{T_0}^T
 \langle Z^{I^\prime } \varPsi _{Tx}(\tau ) , \frac {\sin (T-\tau )B}B
 \psi
\rangle Z^{I^{\prime\prime } }  (\varSigma _T -\frac \varSigma \tau
) \frac {d\tau }\tau \tag 1$$ with $|I^\prime |+|  I^{\prime\prime
}|=|  I|$. We bound pointwise
$$|(1)|\le C    \int _{T_0}^T \| Z^{I^\prime } \varPsi _{Tx}(\tau )
\| _{L^1_x} (T-\tau )^{-\frac 12} |  Z^{I^{\prime\prime } } \frac
1\tau
 (\varSigma _T -\frac \varSigma \tau )| d\tau .$$ Since
$ |I^{\prime  }|+|I^{\prime\prime }|\le N^\prime $,     by Lemma 8.6
and by (6.1)  we have
$$\aligned & \| (1)\| _{L^2_Y} \le C   \int _{T_0}^T  (\|  \varPsi _{Tx}(\tau ) \| _{H^{N^\prime
}_YL^2_x} + \frac {C\epsilon }{\tau ^2})
 (T-\tau )^{-\frac 12}  \tau ^{-\frac 12}
(\| \partial _T \varSigma _{N^\prime }\| _2+\frac 1\tau
\| \varSigma _{N^\prime }\| _2)d\tau.
\endaligned$$
By (6.1) and Lemma 11.3 we conclude   that
$$  \| (2)\| _{L^2_Y} \le \frac{C(\mu ^\prime )}{\delta } \epsilon ^2 T^{-\frac 12+\delta }.$$
 In this way we
conclude Lemma 12.5. \qed

\bigskip

\head \S 13 Closure of the inequalities. \endhead \noindent We
conclude that (6.2) is a consequence of Lemmas 10.2, 11.1, 11.3,
12.2 and 12.4.
\bigskip

\head \S 14 Proofs of the Lemmas   9.2 and 9.3. \endhead

Recall $$\aligned &R_\varPsi =  \Delta_{hyp}\varPsi +T\Cal
Q_0(\varSigma ,   \varSigma  ) P_c \th ^{\prime \prime }  -
 2T \Cal Q_0(\varSigma ,\Cal A) P_c\phi  ^\prime  +\Cal A\Cal
Q_0( \varSigma ,\varSigma )P_c  \phi   ^{\prime \prime }+\\& +P_c
F_3(\varPsi,\Cal{A})+ Q_0( \varSigma ,\varSigma )  P_c\varPsi _{xx}
+ T P_c(  \varPsi_{x}+ \Cal A
 \, \phi  ^\prime )P\varSigma  -T^2
\sum _{\alpha \neq T}  r_{x\alpha } \varPsi_{x\alpha } .
\endaligned
$$
By (4) in Lemma 8.4 we have
$$
\| Z^IB^{m+1}T \Cal Q_0(\varSigma ,\Cal A) P_c\phi  ^\prime \|
_{L^2_{xY}\cap L^2_YL_x^{\frac{p}{p-1}}}\le  C(\mu ^\prime )
\epsilon ^2
 T^{\frac{3}{4}}.
$$

 We will  prove for $|I|+m\le N^\prime $
$$\| Z^IB^{m+1}\left [ R_{\varPsi} +2T \Cal Q_0(\varSigma ,\Cal A) P_c\phi  ^\prime \right ]\| _2\le \left ( C(\mu ^\prime ) \epsilon +C \right )\epsilon
 T^{\frac{3}{4}-\frac{1}{2}+\frac{1}{p}} \tag 14.1$$ which is
stronger than Lemma 9.2. $\| Z^I\Delta _{hyp}\varPsi  \| _2  \le C
T^{C(\mu ^\prime ) \epsilon }\epsilon
   $ by (9.3). For the other terms in the first line, except the   already
discussed
 $-2T \Cal Q_0(\varSigma ,\Cal A) P_c\phi
^\prime $, we have:$\|  Z^I T\Cal Q_0( \varSigma , \varSigma )\| _ 2
\le   C(\mu ^\prime ) \epsilon ^2 T^{ 2\delta }   $ by (6) Lemma
8.4; $\|  Z^I [\Cal A \Cal Q_0( \varSigma , \varSigma )]\| _ 2 \le
C(\mu ^\prime ) \epsilon ^2 T^{ 2\delta -1}   $ by (8) Lemma 8.4. By
(1) in Lemma 7.4 and by Lemma 8.1  we have $ \| Z^I B^m F_j
(\varPsi, \Cal A) \| _2\le {C( \mu^\prime)\epsilon ^2} T^{{C(
\mu^\prime)\epsilon}} $. We have $\| Z^I B^m Q_0( \varSigma
,\varSigma ) P_c\varPsi _{xx}\| _2\le {C( \mu^\prime)\epsilon ^3}
T^{ 2\delta + C( \mu^\prime)\epsilon -1 }$ by (6) in Lemma 8.4 and
Lemma 8.1. By Lemmas 8.1 and 8.5 we have  $\| Z^I B^mT P_c(
\varPsi_{x}+ \Cal A
 \, \phi  ^\prime )P\varSigma \| _2\le {C( \mu^\prime)\epsilon ^3}
T^{ 2\delta + C( \mu^\prime)\epsilon -\frac{1}{2} }$. Next, recall
$$\aligned & r_{xx} =\frac 1{T^2} \Cal Q_0( \varSigma ,
\varSigma ) \quad ,
  \quad
r_{xR}=r_{Rx} =\frac 2{T^3}  \varSigma _R \, ,   r_{xY_j} =r_{Y_jx}=
\frac {2\varSigma _\theta }{T^3 \sinh (R)} a_j(Y)
\endaligned $$
   where
$a_1(Y)   =    -\frac{Y_2}{ \sinh  (R)} $ and $a_2(Y)   =\frac{Y_1}{
\sinh  (R)}  .$ Then: by (6) Lemma 8.4 and by Lemma 8.1 we have $\|
Z^I B^m r_{xx}\varPsi_{xx}\| _2\le {C( \mu^\prime)\epsilon ^3} T^{
2\delta + C( \mu^\prime)\epsilon -1 }$; by (6.1) and by Lemma 8.1 we
have in the other cases $\| Z^I B^m r_{x\alpha }\varPsi_{x\alpha }\|
_2\le {C( \mu^\prime)\epsilon ^2} T^{  \delta - \frac{1}{2}+ C(
\mu^\prime)\epsilon   }.$

Turning to Lemma 9.3, we have proved (14.1) for each single term in
the formula for $R_{\varPsi}$ except for $2T \Cal Q_0(\varSigma
,\Cal A) P_c\phi  ^\prime$. For each term with a cutoff in $x$,  the
estimate translates automatically in the estimate required for Lemma
9.3. For each term linear in $\varPsi$ we can use Lemma 8.6. For the
remaining term, by (2) in  Lemma 7.4 we get
$$ \|
Z^I B ^{m+1}  (\varPsi  ^3) \| _{L^2_Y L^{\frac{p}{p-1}} _x}  \le C(
\mu^\prime) \epsilon \left (  \epsilon + T^{-\frac{1}{2}+2\delta }\|
(\varPsi , \Cal A)_N \| _2\right )
 \| (\varPsi , \Cal A)_N \| _2   .
  $$
  The desired bound follows from Lemma 8.1.

\bigskip
ACKNOWLEDGMENTS. I wish to thank Professor Fabio Zanolin for his
hospitality at the University of Udine.

\enddocument